\documentclass[12pt]{article}

\usepackage{amsmath}
\usepackage{amssymb}

\newtheorem{Th}{Theorem}[section]
\newtheorem{Prop}{Proposition}[section]

\newtheorem{Lemma}{Lemma}[section]
\newtheorem{Coro}{Corollary}[section]

\newtheorem{Rem}{Remark}[section]

\newcommand{\finishproof}{\hfill $\Box$ \vspace{3mm}}
\newcommand{\R}{\mathbb{R}}
\newcommand{\Z}{\mathbb{Z}}
\newcommand{\C}{\mathbb{C}}

\newcommand{\W}{{\mathcal W}}

\newcommand{\spec}{\mathop{\tt spec}}

\newcommand{\la}{\lambda}

\newcommand{\om}{\omega}
\newcommand{\Om}{\Omega}
\newcommand{\Ga}{\Gamma}
\newcommand{\De}{\Delta}

\begin{document}

\title{On normalized differentials on hyperelliptic curves of infinite genus}

\author{T. Kappeler and P. Topalov}

\date{}

\maketitle

\abstract{We develop a new approach for constructing normalized differentials
on hyperelliptic curves of infinite genus and
obtain uniform asymptotic estimates for the distribution of their zeros.}


\section{Introduction}\label{1. Introduction}

Much of the analysis of closed Riemann surfaces is based on the Riemann bilinear
relation. Given a canonical basis $A_1, B_1, \ldots , A_g, B_g$, of the homology group
$H_1(X,{\mathbb Z})$ of a closed Riemann surface $X$ of genus $g$, it reads
   \[ \int_X \omega \wedge \eta = \sum ^g_{m=1} \Big( \int _{A_m} \omega
      \int _{B_m} \eta - \int _{A_m} \eta \int _{B_m} \omega \Big)
   \]
where $\omega$ and $\eta $ are arbitrarily chosen smooth closed $1$-forms on $X$. 
As a consequence one obtains the following vanishing theorem: for any holomorphic
$1$-form $\omega $ on $X$ with vanishing $A$-periods, 
$\int _{A_m} \omega =0$ $\forall 1 \leq m \leq g$, one has 
$-\frac{1}{2 i}\int_X \omega \wedge \overline{\omega } = 0$, and hence $\omega = 0$.
Furthermore it follows from Hodge theory that the space of holomorphic differentials
on $X$ is a complex vector space of dimension $g$ admitting a basis
$\omega _1, \ldots , \omega _g$ such that $\int _{A_m} \omega _n = \delta _{mn}$
for any $1 \leq m, n \leq g$. By the vanishing theorem, such a basis is unique.
The period matrix $B_X = \big(\int_{B_m}\omega_n\big)_{1 \leq m, n \leq g}$, 
which enters the definition of the theta function associated to $X$,
is known to be symmetric and has the property that $\mathop{\rm Im} B_X $ is sign definite.

However, for many applications to integrable PDEs one needs to consider {\em open}
Riemann surfaces of {\em infinite genus}, a subject pioneered by Ahlfors and Nevanlinna -- see the
monographs \cite{AS} and \cite{FKT} as well as references therein.
Unfortunately, it is not sufficiently developed for our purposes. In
particular for applications to the focusing nonlinear Schr\"odinger (NLS) equation
we need to establish a vanishing theorem for holomorphic $1$-forms on two
sheeted open Riemann surfaces of infinite genus which are not necessarily
$L^2$-integrable. In view of these applications we formulate our results for the specific
Riemann surfaces involved. However, our method is quite general and can be directly
applied for studying Riemann surfaces related to other non-linear equations.

Consider the NLS system of equations in one space dimension with 
periodic boundary conditions,
\begin{equation}
\left\{
\begin{array}{c}
i \partial_t \varphi_1 = - \partial^2_x \varphi _1 + 2 \varphi^2_1\varphi _2,\\
-i \partial_t \varphi_2 = - \partial^2_x \varphi_2 + 2 \varphi ^2_2 \varphi_1,
\end{array}
\right.
\end{equation}
where $\varphi = (\varphi _1, \varphi _2)$ is in $L^2_c:= L^2({\mathbb T},
{\mathbb C}) \times L^2({\mathbb T},{\mathbb C})$ and 
${\mathbb T} ={\mathbb R} / {\mathbb Z}$. If $\varphi _2 = \overline{\varphi _1}$
[$\varphi _2= - \overline{\varphi _1}$], the system is referred to as defocusing [focusing]
NLS equation. By Zakharov and Shabat \cite{ZS} the NLS system admits a Lax
pair formulation, $\partial _t L= [B, L]$ where $L(\varphi )$ denotes the
Zakharov-Shabat (ZS) operator
   \[ L(\varphi ) := i \begin{pmatrix} 1 & 0 \\ 0 &-1 \end{pmatrix} \partial_x
      + \begin{pmatrix} 0 &\varphi_1 \\ \varphi_2 & 0 \end{pmatrix}\,.
   \]
Associated to this operator is the curve
   \[ {\mathcal C}_\varphi := \big\{ (\lambda,w) \in {\mathbb C}^2\,\big|\, w^2 =
      \Delta ^2 (\lambda,\varphi) - 4 \big\}
   \]
where $\Delta(\lambda,\varphi)$ is the discriminant of $L(\varphi)$
(cf. Section~\ref{2. Preliminaries}). It is known (see e.g. \cite{GKP})
that for any given $\varphi\in L^2_c$, the entire function
$\Delta^2(\lambda,\varphi)-4$ vanishes at $\lambda\in{\mathbb C}$ iff 
$\lambda$ is a periodic eigenvalue of $L(\varphi)$, i.e., an eigenvalue of $L(\varphi)$,
considered on the interval $[0,2]$ with periodic boundary conditions. 
In addition, the algebraic multiplicity of $\lambda$ as a root of 
$\Delta(\lambda,\varphi)^2-4$ coincides with the algebraic multiplicity of it as
a periodic eigenvalue of $L(\varphi)$ (see \cite{KLT2}).
Note that $L(\varphi )$ has a compact resolvent
and hence its spectrum is discrete. In this paper we consider potentials
$\varphi$ in $L^2_c$ so that each periodic eigenvalue $\lambda$
of $L(\varphi)$ has algebraic multiplicity $m(\lambda )$ at most two.
Denote the subset of such elements by $L^2_\bullet $. Then $L^2_\bullet$ is
open, dense, and contains the zero potential (\cite{KLT2}). 
In the case of the defocusing NLS equation, $L(\varphi)$ is self-adjoint. 
Then the algebraic multiplicity of an eigenvalue coincides with the geometric one
which is at most two, and hence $\varphi\in L^2_\bullet$ in the defocusing case.
For $\varphi \in L^2_\bullet$, the periodic eigenvalues of $L(\varphi )$ can be listed
as a sequence of distinct pairs, $\lambda ^-_k(\varphi )$, $\lambda ^+_k(\varphi )$, 
$k \in{\mathbb Z}$, so that $\lambda^\pm_k(\varphi) = k\pi + \ell ^2(k)$,
i.e., $\sum_{k \in {\mathbb Z}} |\lambda ^\pm _k(\varphi ) - k\pi |^2< \infty $,
and $\lambda_k^-=\lambda_k^+$ iff $\lambda_k^-$ has algebraic multiplicity two,
$m(\lambda)=2$ (Section~\ref{2. Preliminaries}).

Let $Z_\varphi $ denote the subset of periodic double eigenvalues of
$L(\varphi )$,
   \[ Z_\varphi := \big\{ \lambda \in \mbox{spec } L(\varphi )\,\big|\, m(\lambda ) =
      2 \big\} .
   \]
Then 
\[
{\mathcal C}^\bullet_\varphi := {\mathcal C} _\varphi \setminus 
\big\{(\lambda,0)\,\big|\, \lambda \in Z_\varphi\big\}
\] 
is a two sheeted open Riemann surface. Generically it is a surface of infinite genus 
(\cite{KLT2}). 
Our aim is to prove a vanishing theorem for holomorphic differentials on 
${\mathcal C}^\bullet_\varphi $ which are not necessarily $L^2$-integrable and to
construct a family of normalized holomorphic differentials $\omega _n$, $n \in {\mathbb Z}$,
on ${\mathcal C}^\bullet_\varphi$ with respect to an appropriately chosen infinite set of
cycles on ${\mathcal C}^\bullet_\varphi $.
In addition we want to get asymptotic estimates of the zeros of these differentials. 
The cycles are defined as follows: for any given potential
$\psi\in L^2_\bullet\,$ list its periodic eigenvalues in pairs, 
$\lambda ^-_k,\lambda^+_k$, $k \in {\mathbb Z}$, as discussed above. 
It is shown in Section~\ref{2. Preliminaries} that there exist an open neighborhood
${\mathcal W}$ of $\psi$ in $L^2_\bullet$ and a family of simple, closed, smooth,
counterclockwise oriented curves $\Gamma _m, m \in {\mathbb Z}$, so that
the closures of the domains in ${\mathbb C}$, bounded by the $\Gamma_m$
are pairwise disjoint and for any $\varphi \in {\mathcal W}$ and 
$m \in{\mathbb Z}$ the open domain bounded by $\Gamma_m$ contains the pair 
$\lambda^-_m(\varphi )$, $\lambda^+_m(\varphi )$ but no other periodic eigenvalues
of $L(\varphi)$. Denote by $A_m$ the cycle on the canonical sheet
${\mathcal C}^c_\varphi $ of ${\mathcal C}_\varphi$ 
(cf. Section~\ref{2. Preliminaries}) so that $\pi (A_m) = \Gamma _m$ where
$\pi : {\mathcal C}_\varphi \rightarrow {\mathbb C}$, 
$(\lambda,w)\mapsto\lambda $. We are then looking for holomorphic differentials
$\omega _n$ on ${\mathcal C}^\bullet _\varphi $ so that 
$\int _{ A_m} \omega _n= 2\pi \delta _{mn}$ for any $m,n\in {\mathbb Z}$. 
In addition we want to prove a vanishing theorem for holomorphic differentials on 
${\mathcal C}^\bullet_\varphi $ with vanishing $A$-periods which are not necessarily
$L^2$-integrable.
For an arbitrary entire function $\zeta : {\mathbb C} \rightarrow {\mathbb C}$
with $\zeta \big\arrowvert _{Z_\varphi } = 0$, let 
\[
\omega_\zeta :=\frac{\zeta (\lambda )}{\sqrt{\Delta ^2(\lambda , \varphi ) - 4}}\,d\lambda\,.
\]
Then $\omega _\zeta $ is a holomorphic $1$-form on ${\mathcal C}^\bullet_\varphi$ 
which is locally square integrable. More precisely, for any $r> 0$,
   \[ V(r):= - \frac{1}{2i} \int _{X_r} \omega _\zeta \wedge \overline{\omega
      _\zeta } < \infty
   \]
where 
\[
X_r:= \pi ^{-1} \big( \{ \lambda \in {\mathbb C}\,\big|\, |\lambda | \leq r \}\big) 
\cap{\mathcal C}^\bullet _\varphi
\] 
and where the orientation is induced by the complex structure on ${\mathbb C}$. 
Using polar coordinates $(\rho,\theta)$ on ${\mathbb C}$ and taking into account
that ${\mathcal C}_\varphi$ is a two sheeted curve one has
   \[ V(r) = 2 \int ^r_0 \int ^{2\pi }_0 \Big|\frac{\zeta (\rho
      e^{i\theta })}{\sqrt{\Delta ^2(\rho e^{i \theta }) - 4}} \Big|^2 
    \rho\,d\theta d \rho\,.
   \]
In particular, if $\zeta \not\equiv 0$ then $V(r) > 0$ for any $r > 0$ and
$V(r)$ is strictly increasing. Note that $\omega _\zeta $ might not be
$L^2$-integrable as it could happen that 
$\lim\limits_{r \rightarrow \infty } V(r) = \infty $.

\begin{Th}\label{Theorem 1.0} 
Let $\varphi\in\W$ with $\W$ and $(A_m)_{m \in{\mathbb Z}}$
given as above and let $\zeta : {\mathbb C} \rightarrow {\mathbb C}$ be
an entire function with $\zeta \not\equiv 0$ and $\zeta \big\arrowvert _{Z_\varphi } = 0$. 
If $\int _{A_m} \omega _\zeta = 0$ for any $m \in {\mathbb Z}$, then there exists $C > 0$
so that 
\[
V(r) \geq C\,r^{2/\pi }
\] 
for any $r \geq 1$.
\end{Th}

Theorem~\ref{Theorem 1.0} leads to the following vanishing theorem.

\begin{Th}\label{Theorem 1.1}
Let $\varphi\in\W$ with $\W$ and $(A_m)_{m \in {\mathbb Z}}$
given as above and let $\zeta : {\mathbb C} \rightarrow {\mathbb C}$ be entire
with $\zeta \big\arrowvert _{Z_\varphi } = 0$. If $\int _{A_m} \omega _\zeta= 0$
for any $m \in {\mathbb Z}$ and 
\[
V(m\pi) = o(m^{2/\pi })\,\,\,\mbox{as}\,\,\,m\rightarrow \infty 
\]
then $\zeta\equiv 0$ and hence $\omega_\zeta\equiv 0$.
\end{Th}

\begin{Rem}\label{Remark 1}
The conclusion of Theorem~\ref{Theorem 1.1} no longer holds
when the assumption $V(m\pi ) = o(m^{2/\pi })$ is dropped -- 
see the form $\hat\om$ defined in \eqref{eq:hat_omega}.
\end{Rem}

\begin{Rem}\label{Remark 2}
As $V(r)$ is increasing, the conditions $V(m\pi) = o\big( m^{2/\pi}\big)$ as 
$m \rightarrow \infty $ is equivalent to $V(r) = o \big( r^{2/\pi })$ as 
$r \rightarrow \infty $.
\end{Rem}

Next we state our result on normalized differentials on ${\mathcal C}^\bullet_\varphi$ 
and describe features of them needed in our studies of the focusing NLS equation.

\begin{Th}\label{Theorem 1.2} 
For any $\psi\in L^2_\bullet$ there exist an open neighborhood ${\mathcal W}$ of $\psi$ in $L^2_\bullet$, 
cycles $A_m$, $m \in {\mathbb Z}$, as above and analytic functions 
$\zeta_n : {\mathbb C} \times{\mathcal W} \rightarrow {\mathbb C}$, 
$n \in {\mathbb Z}$, so that for any $\varphi \in {\mathcal W}$ and $n \in {\mathbb Z}$, 
the holomorphic differential 
$\omega _n = \frac{\zeta_n(\lambda,\varphi)}{\sqrt{\Delta^2(\lambda,\varphi)-4}}\,d\lambda$
on ${\mathcal C}^\bullet_\varphi$
satisfies
   \[ 
       \frac{1}{2\pi } \int _{A_m} \omega _n = \delta _{nm} \quad \forall m\in {\mathbb Z}\,.
   \]
In addition, there exists $N\ge 1$ such that for any $n\in\Z$ and for any $\varphi\in\W$
the entire function $\zeta_n(\cdot,\varphi)$ admits infinitely many zeros. When listed appropriately, 
these zeros $\sigma ^n_k$, $k \in {\mathbb Z} \backslash \{ n \} $, satisfy:
\begin{itemize}
\item[(i)] For any $|k| > N$, $k \neq n$, $\sigma^n_k$ is the only zero of
$\zeta_n(\cdot,\varphi)$ in the disk $D_k(\pi/4)$ and the map 
$\sigma^n_k : {\mathcal W}\to D_k(\pi / 4)$ is analytic.
Furthermore, for any $|k| \leq N$, $k \neq n$, $\sigma ^n_k \in D_0(N\pi + \pi /4)$.
There are {\em no} other zeros of $\zeta_n(\cdot,\varphi)$ in $\C$.

\item[(ii)] For any $|k| > N$, $k\ne n$,
\begin{equation}\label{eq:sigma-asymptotics}
\sigma^n_k=\frac{\lambda^-_k+\lambda^+_k}{2}+\big(\lambda^+_k-\lambda^-_k\big)^2\,\ell^2(k)
\end{equation}
uniformly in $n\in\Z$ and locally uniformly in $\W$.
\end{itemize}
Moreover, $\zeta_n(\cdot,\varphi)$ admits the product representation
   \[ \zeta _n(\lambda , \varphi ) = - \frac{2}{\pi _n} \prod _{k \not= n}
      \frac{\sigma ^n_k - \lambda }{\pi _k}
   \]
where 
   \[ \prod _{k \neq n}\frac{\sigma^n_k-\lambda}{\pi _k}:=
      \lim_{K\to\infty}\prod_{|k|\leq K, k\neq n} 
          \frac{\sigma^n_k-\lambda}{\pi_k}
   \]
and
\[
\pi_k:=\left\{
\begin{array}{l}
k\pi,\,\,\,\,k\ne 0\\
1,\,\,\,\,\,\,\,k=0
\end{array}
\right..
\]
For any $\lambda^\pm_k\in Z_\varphi$ with $k\neq n$, $\zeta_n(\lambda^\pm_k,\varphi)=0$, 
and $\omega_n$ is $L^2$-integrable in sufficiently small punctured neighborhoods in 
${\mathcal C}_\varphi^\bullet$ of the point $(\lambda^\pm_k,0)$. If, however 
$\lambda^\pm_n\in Z_\varphi $, then $\omega _n$ is not $L^2$-integrable in any punctured
neighborhood in ${\mathcal C}_\varphi^\bullet$ of the point $(\lambda^\pm_n,0)$.
\end{Th}

\begin{Rem} Using the product representation for $\zeta _n(\lambda)$, the
asymptotic estimates for the $\sigma ^n_k$'s, and estimates on infinite
products in \cite[Lemma C.5]{GKP} one can show that for any $n\in{\mathbb Z}$ and 
$\varphi \in {\mathcal W}$ there exists $C > 0$ so that
$-\frac{1}{2i}\int_{X_r\setminus X_{|n|\pi+\pi/4}}\om_n\wedge\overline{\om_n}\ge C \log r$
for any $r\ge |n|\pi+\pi/4$. Hence $\omega _n$ is never $L^2$-integrable.
\end{Rem}

\begin{Rem} The uniformity statement in the asymptotic formula \eqref{eq:sigma-asymptotics}
means that for any $\varphi \in {\mathcal W}$ there is a neighborhood ${\mathcal V}$ of
$\varphi$ in $\W$ and a constant $C > 0$ so that for any $n\in{\mathbb Z}$ and for any
$\nu\in{\mathcal V}$ there are constants $(c^n_k)_{|k|>N, k\ne n}$, $c^n_k\ge 0$, so that
for any $|k|>N$, $k\ne n$,
   \[ 
        \left|\sigma ^n_k - \frac{\lambda ^+_k + \lambda ^-_k}{2}\right|
         \leq |\lambda ^-_k - \lambda ^-_k |^2 c^n_k
   \] 
and $\sum\limits_{k \not= n} |c^n_k|^2 \leq C$.
\end{Rem}
Let $\W$ and $A_m$, $m\in\Z$, be as above. Consider the holomorphic 1-form on 
${\mathcal C}^\bullet_\varphi$,
\[
\Om^*:=\frac{\dot\Delta(\lambda)}{\sqrt{\Delta^2(\lambda)-4}}\,d\lambda,
\]
where $\dot\Delta(\lambda):=\partial_\lambda\Delta(\lambda,\varphi)$.
By Lemma \ref{lem:zero_form} in Section \ref{2. Preliminaries} below one can choose $\W$ and $N\ge 1$ in
Theorem \ref{Theorem 1.2} so that for any $|m|\ge N$,
\begin{equation}\label{eq:zero-cycles}
\int_{A_m}\Om^*=0\,. 
\end{equation}
Hence,
\begin{equation}\label{eq:hat_omega}
{\hat\om}:=\Om^*-\frac{1}{2\pi}\sum_{|m|<N}a_m\om_m
\end{equation}
with $a_m:=\int_{A_m}\Om^*$ and $\om_m$ given by Theorem \ref{Theorem 1.2},
is a holomorphic 1-form on ${\mathcal C}^\bullet_\varphi$ satisfying $\int_{A_m}{\hat\om}=0$
for any $m\in\Z$. In particular, the conclusion of Theorem \ref{Theorem 1.1} does {\em not} hold if
the growth condition is dropped.

Clearly, the differentials $\omega _n$ of Theorem~\ref{Theorem 1.2} are
unique within the class of holomorphic differentials obtained by perturbations
of the type defined by Theorem~\ref{Theorem 1.1}. However, without some
conditions on the behaviour of the differentials near infinity, one cannot
expect uniqueness. Indeed, any sequence of holomorphic differentials of the
form ${\tilde\omega}_n:=\omega _n + c_n \hat \omega$, $c_n \in {\mathbb C}$, satisfies
$\frac{1}{2\pi } \int _{A_m}{\tilde\omega}_n = \delta _{mn}$ for any $m, n \in {\mathbb Z}$.

Besides Theorem~\ref{Theorem 1.1}, a key ingredient in the proof of
Theorem~\ref{Theorem 1.2} is a novel ansatz for the entire functions $\zeta_n$, $n\in\Z$, 
leading to a {\it linear} system of equations. A detailed outline of
the proof of Theorem~\ref{Theorem 1.2} is given in Section~\ref{3. Outline of proof of Theorem 1.2}.

\hspace{0.5cm}

\noindent{\em Applications:} In \cite{KLT}, Theorem 1.2 is used to construct
action-angle variables for the focusing NLS equation, significantly extending
previous our results obtained in \cite{KLTZ} near the zero potential. 
See \cite{AbMa} for related results for $1$-gap and $2$-gap potentials and
\cite{DMN}, \cite{IK}, \cite{PT}, \cite{VN} for finite gap potentials. Such
coordinates allow to obtain various results concerning well-posedness for
these equations and study their (Hamiltonian) perturbations -- see e.g.
\cite{KP}, \cite{KT2}, where results in this
direction have been obtained for the KdV equation.

\hspace{0.5cm}

\noindent{\em Related results:} (1) In \cite{GKP} (cf. also \cite{MV}) the
case of the defocusing NLS equation is treated, i.e., the case where 
$\varphi = (\varphi_1, \varphi_2)$ satisfies $\varphi_2 = \overline{\varphi_1}$. 
In this case $\varphi \in L^2_\bullet $ and $L(\varphi )$ is self-adjoint, 
hence its periodic spectrum is real. More precisely the eigenvalues can be 
listed in such a way that
   \[ \cdots < \lambda ^-_k \leq \lambda ^+_k < \lambda ^-_{k + 1}
      \leq \lambda ^+_{k + 1} < \cdots \quad \mbox { and } \quad \lambda ^\pm 
      _k = k\pi + \ell ^2(k) .
   \]
It then follows that zeros of $\zeta _n$ are confined to the closed gaps
$[\lambda ^-_k, \lambda ^+_k]$ with $k \not= n$. More precisely, as
$\Delta (\lambda )$ is real valued and $|\Delta (\lambda )| > 2$ on the
open gaps $(\lambda ^-_k, \lambda ^+_k), k \in {\mathbb Z}$, one easily
deduces that for any $n \in {\mathbb Z}, \, \zeta _n(\lambda )$ must have a
zero $\sigma ^n_k$ in any closed gap $[\lambda ^-_k, \lambda ^+_k]$ with
$k \not= n$. Using the implicit function theorem it is shown in
\cite{GKP} that normalized differentials can be constructed in a sufficiently small
neighborhood $W \subseteq L^2_\bullet $ of such a $\varphi $. 
Note that if $\varphi \in L^2_\bullet $ is arbitrary we have no a priori
knowledge on the zeros of $\zeta _k$ and hence one {\em cannot} apply the
implicit function theorem approach in \cite{GKP} as well as the method in \cite{MV}.

In \cite{KLT1} we consider finite gap potentials $\varphi \in L^2_\bullet$ so that
the ZS operator $L(\varphi)$ is {\em not} necessarily self-adjoint. Listing the
periodic eigenvalues of $L(\varphi )$ in pairs $\lambda ^-_k(\varphi ), \lambda ^+_k(\varphi )$, 
$k \in {\mathbb Z}$, as above it means that the subset $J_\varphi $ consisting of all
$k \in {\mathbb Z}$ with $\lambda ^-_k(\varphi ) \not= \lambda ^+_k(\varphi )$ is finite.
Using again the implicit function theorem it is shown in \cite{KLT1}
that normalized differentials can be constructed in a sufficiently small neighborhood
of such a finite gap potential. Unfortunately, the method in \cite{KLT1} does {\em not} work
if $\varphi$ is an arbitrary not necessarily finite gap potential in $L^2_\bullet$.

\vspace{0.2cm}

(2) In the case of the KdV equation on the circle, the relevant operator in
the Lax pair is the Hill operator. For potentials with sufficiently
small imaginary part, normalized differentials have been constructed e.g. in
\cite{KP}, using the implicit function theorem and the same method as in \cite{GKP}. 
In the case where the Hill operator has simple periodic spectrum, the corresponding spectral curve
is a Riemann surface of infinite genus and the existence of normalized holomorphic differentials
can be proved by Hodge theory using the fact that in this case these
differentials are $L^2$-integrable -- see \cite{FKT}, \cite{MT}.
Note however that these arguments do {\em not} provide the asymptotic estimates
of the zeros nor the analytic dependence on the potential. But even for
the existence part this approach would not work for the ZS operator as
the differentials of Theorem~\ref{Theorem 1.2} are never $L^2$-integrable
on ${\mathcal C}^\bullet _\varphi $.

\hspace{0.5cm}

\noindent{\em Organization:} In the preliminary Section~\ref{2. Preliminaries}
we introduce additional notation and review the spectral properties of ZS
operators needed throughout the paper. In Section~\ref{4. Proof of Theorem 1.0 
and Theorem 1.1}
we prove Theorem~\ref{Theorem 1.0} and Theorem~\ref{Theorem 1.1}
whereas in Section~\ref{3. Outline of proof of 
Theorem 1.2} we give an outline of the
proof of Theorem~\ref{Theorem 1.2}. Its details are then presented in the
remaining sections.


\section{Preliminaries}\label{2. Preliminaries}
In this section we introduce some more notation and review properties of the
Zakharov-Shabat operator $L(\varphi )$, introduced in Section~\ref{1. Introduction}.
For $\varphi \in L^2_c$ and $\lambda \in {\mathbb C}$, let 
$M(x) \equiv M(x, \lambda , \varphi )$ denote the fundamental $2 \times 2$ matrix of 
$L(\varphi ), \, L(\varphi )M(x) = \lambda M(x)$, satisfying the initial condition
$M(0, \lambda , \varphi ) = Id_{2 \times 2}$. Let $\R_{\ge 0}:=\{x\in\R\,|\,x\ge 0\}$.
It is well known that 
$M : {\mathbb R} _{\geq 0} \times {\mathbb C} \times L^2_c \rightarrow {\mathbb C}^{2 \times 2}$ 
is continuous and for any $x$ fixed, 
$M(x, \cdot , \cdot ) : {\mathbb C} \times L^2_c \rightarrow {\mathbb C}^{2 \times 2}$ is analytic --
see e.g. \cite{GKP}, Chapter I.

\hspace{0.5cm}

\noindent {\em Periodic spectrum:} Denote by $\spec L(\varphi )$ the spectrum
of $L(\varphi )$ with domain
   \[ \mbox{dom}_{per}L(\varphi ):= \{ F \in H^1_{loc}\times H^1_{loc}: F(2) =
      F(0) \}
   \]
where $H^1_{loc} \equiv H^1_{loc}({\mathbb R}, {\mathbb C})$. As
$L(\varphi )$ has a compact resolvent, the periodic spectrum of $L(\varphi )$
is discrete. It has been analyzed in great detail.

The discriminant $\Delta (\lambda ) \equiv \Delta (\lambda , \varphi )$ of
$L(\varphi )$ is defined to be the trace of $M(1, \lambda , \varphi )$, 
$\Delta(\lambda ) = tr M(1, \lambda )$. It is straightforward to see that $\lambda $ is
a periodic eigenvalue of $L(\varphi )$ iff $\Delta ^2(\lambda , \varphi ) - 4 =
0$. Clearly, $\Delta : {\mathbb C} \times L^2_c \rightarrow {\mathbb C}$ is
analytic.

We say that $a, b \in {\mathbb C}$ are lexicographically ordered, $a \preccurlyeq b$, iff 
$[\mathop{\rm Re}(a)< \mathop{\rm Re}(b)]$ or $[\mathop{\rm Re}(a) = \mathop{\rm Re}(b)$ and 
$\mathop{\rm Im}(a) \leq \mathop{\rm Im}(b)]$. Similarly, $a \prec b$
iff $a \preccurlyeq b$ and $a \not= b$. The following two propositions are
well known -- see e.g. \cite[Section 3]{GKP} or references therein. For any
$k \in {\mathbb Z}$ and $r > 0$ denote by $D_k(r)$ the disk
   \[ D_k(r):= \{ \lambda \in {\mathbb C}\,|\, |\lambda - k\pi | < r\} .
   \]

\begin{Prop}\label{Proposition 2.1} 
For any $\psi\in L^2_c$ there exist an open neighborhood ${\mathcal W}$ of $\psi$ in $L^2_c$ and 
an integer $N_0 \geq 1$ such that for any $\varphi \in {\mathcal W}$ the following holds:
\begin{itemize}
\item[(i)] For any $k \in {\mathbb Z}$ with $|k| \geq N_0$, the disk $D_k(\pi/6)$ contains precisely
two periodic eigenvalues $\lambda^-_k(\varphi ) \preccurlyeq \lambda ^+_k(\varphi )$ of $L(\varphi )$
and one zero $\dot \lambda _k(\varphi )$ of 
$\dot \Delta (\lambda , \varphi )= \partial_\lambda \Delta (\lambda , \varphi )$ 
(all counted with their algebraic multiplicities).
\item[(ii)] The disk $D_0\big((N_0 - 3/4)\pi\big)$ contains precisely $4N_0 - 2$
periodic eigenvalues of $L(\varphi )$ and $2N_0 - 1$ zeros of $\dot \Delta
(\lambda , \varphi )$ (all counted with their algebraic multiplicities).
\item[(iii)] There are no other periodic eigenvalues of $L(\varphi )$ and
no other zeros of $\dot \Delta (\lambda , \varphi )$ than the ones listed
in (i) and (ii).
\end{itemize}
\end{Prop}

\begin{Prop}\label{Proposition 2.2}
Let ${\mathcal W} \subseteq L^2_c$ be given by
Proposition~\ref{Proposition 2.1}. For any $\varphi \in {\mathcal W}$, the
periodic eigenvalues $(\lambda ^\pm _k)_{|k| \geq N_0}$ and the zeros
$(\dot \lambda _k)_{|k| \geq N_0}$ satisfy the asymptotic estimates
   \[ \lambda ^\pm _k = k\pi + \ell ^2(k) \mbox { and } \dot \lambda _k
      = k\pi + \ell ^2(k)
   \]
locally uniformly in $\W$. More precisely, it means that
$\sum _{|k| \geq N_0} |\dot \lambda _k - k\pi |^2$ is locally
bounded in $\W$. Similar statements hold for $\lambda^\pm_k$.
\end{Prop}

Take  $\psi\in L^2_\bullet$ and construct a neighborhood $\W\subseteq L^2_\bullet$
of $\psi$ in $L^2_\bullet$ so that the statements of Proposition \ref{Proposition 2.1} and
Proposition \ref{Proposition 2.2} hold.\footnote{Recall that $L^2_\bullet$ is open 
and dense in $L^2_c$.}
For any $\varphi\in\W$, in addition to the periodic eigenvalues 
$(\lambda ^\pm _k)_{|k| \geq N_0}$ of the ZS operator $L(\varphi)$ there are $4N_0 - 2$ 
periodic eigenvalues in the disk $D_0\big((N_0 - 3/4)\pi\big)$. We list these eigenvalues in pairs 
$\lambda ^-_k,\lambda ^+_k$, $|k| < N_0$, in an arbitrary way except that 
any double eigenvalue is listed as a pair and for any $|k|<N_0$ the eigenvalues
$\lambda_k^-$ and $\lambda_k^+$ are lexicographically ordered 
$\lambda ^-_k \preccurlyeq \lambda ^+_k$. 
For all integers $|k| < N_0$, choose simple, closed smooth, counterclockwise
oriented curves $\Gamma_k$ contained in the disk $D_0\big((N_0 - 3/4)\pi\big)$
so that the closures of the (open) domains in $\C$, bounded by the
$\Gamma_k$ are pairwise disjoint and for any $|k| < N_0$, the domain bounded
by $\Gamma_k$ contains the pair $\lambda ^\pm_k$, but no other
periodic eigenvalue of $L(\varphi)$. For each $|k| < N_0$ choose a closed smooth
curve $\Gamma_k'$ in the interior of $\Gamma_k$ so that 
$\mbox{dist}(\Gamma_k, \Gamma_k') > 0$ and $\lambda^\pm_k$ are inside $\Gamma'_k$.
By shrinking the neighborhood $\W$ of $\psi$ in $L^2_\bullet$ if necessary,
the choice of $\Gamma'_k$, $|k|<N_0$, can be done uniformly in $\W$, i.e., for any 
$\varphi \in {\mathcal W}$ and for any $|k|<N_0$ the domain bounded by $\Gamma_k'$
contains precisely two periodic eigenvalues of $L(\varphi)$, 
$\lambda ^-_k(\varphi ) \preccurlyeq \lambda ^+_k(\varphi )$. 
Furthermore for any $\varphi \in {\mathcal W}$ and $|k| < N_0$ chose a continuously
differentiable simple curve $G_k\equiv G_k(\varphi)$ inside $\Gamma_k'$ connecting 
$\lambda^-_k(\varphi)$ with $\lambda^+_k(\varphi)$. In the case when
$\lambda^-_k(\varphi)=\lambda^+_k(\varphi)$, $G_k(\varphi)$ is chosen to be the constant curve
$\lambda^-_k(\varphi)$. For $|k|\ge N_0$, we choose $\Gamma_k$ to be the
counterclockwise oriented boundary of the disk $D_k(\pi/4)$
and $G_k$ to be the straight line,
   \[ G_k : [0,1] \rightarrow D_k(\pi/4), \,\, 
            t\mapsto(1 - t)\lambda^-_k(\varphi)+t\lambda^+_k(\varphi)\, .
   \]
Furthermore for $k \in {\mathbb Z}$ and $\varphi \in {\mathcal W}$ we define
\[
\tau_k(\varphi ) := \big(\lambda^-_k(\varphi) + \lambda^+_k(\varphi)\big)/2
\]
and
\[
\gamma _k(\varphi ) := \lambda ^+_k(\varphi )-\lambda ^-_k(\varphi)\,.
\]

\vspace{0.3cm}

\noindent {\em Infinite products:} We say that an infinite product 
$\prod_{k\in {\mathbb Z}} (1 + a_k)$ with $a_k \in {\mathbb C}$, $k \in {\mathbb Z}$,
converges if $\lim_{K \rightarrow \infty } \prod _{|k| \leq K} (1 + a_k)$ exists. 
The limit is then also denoted by $\prod _{k \in {\mathbb Z}} (1 + a_k)$.
The infinite product $\prod_{k \in {\mathbb Z}}(1 + a_k)$
converges absolutely if $\prod_{k \in {\mathbb Z}}(1 + |a_k|)$ converges.

\vspace{0.3cm}

\noindent {\em Product representations:} For any $\varphi \in L^2_c, \Delta
^2 (\lambda , \varphi ) - 4$ and $\dot \Delta (\lambda , \varphi )$ admit
product representations. For any given element in $L^2_\bullet$, choose $N_0$
and ${\mathcal W}$ as in Proposition~\ref{Proposition 2.1}. According to
Proposition~\ref{Proposition 2.1}, for any $\varphi \in {\mathcal W},
\dot \Delta (\lambda , \varphi )$ admits $2N_0 - 1$ zeros in the disk
$D_0\big( (N_0 - \frac{3}{4}) \pi \big) $. For convenience list them in
lexicographic order, $\dot \lambda _{-N_0 + 1}(\varphi ) \preccurlyeq
\dot \lambda _{- N_0 + 2}(\varphi ) \preccurlyeq \ldots \preccurlyeq
\dot \lambda _{N_0 - 1}(\varphi )$. The remaining zeros are listed
as in Proposition~\ref{Proposition 2.1}. The proof of the following
statement can be found in \cite[Lemma 6.5, Lemma 6.8]{GKP}.

\begin{Prop}\label{Proposition 2.3}
For any $\varphi \in {\mathcal W}$ and $\lambda \in{\mathbb C}$
   \[ \Delta ^2(\lambda , \varphi ) - 4 = - 4 \prod _{k \in {\mathbb Z}}
      \frac{(\lambda ^+_k(\varphi ) - \lambda )(\lambda ^-_k(\varphi ) -
      \lambda )}{\pi ^2_k}
   \]
and
   \[ \dot \Delta (\lambda , \varphi ) = 2 \prod _{k \in {\mathbb Z}}
      \frac{\dot \lambda _k(\varphi ) - \lambda }{\pi _k} .
   \]
\end{Prop}

\vspace{0.3cm}

\noindent {\em Standard and canonical roots:} Denote by $\sqrt[+]{z}$
the branch of the square root defined on ${\mathbb C} \backslash
\{ z \in {\mathbb R} : z \leq 0 \}$ by $\sqrt[+]{1} = 1$. For any $a, b \in
{\mathbb C}$, we define the standard root of $(a - \lambda )(b - \lambda )$
by the following relation
   \begin{equation}\label{eq:canonical_root} 
         \sqrt[s]{(a - \lambda )(b - \lambda )} = - \lambda 
             \sqrt[+]{\Big(1 - \frac{a}{\lambda }\Big)\Big(1-\frac{b}{\lambda}\Big)}
   \end{equation}
for all $\lambda \in {\mathbb C} \backslash \{ 0 \}$ such that $\Big\arrowvert
\frac{a}{\lambda } \Big\arrowvert , \Big\arrowvert \frac{b}{\lambda } \Big
\arrowvert \leq 1/2$. Let $G_{[a,b]}$ be a continuous simple curve connecting $a$ and $b$. 
By analytic extension, \eqref{eq:canonical_root} uniquely defines a
holomorphic function on ${\mathbb C} \backslash G_{[a,b]}$ that we call the
{\it standard root} of $(a - \lambda )(b - \lambda )$ on ${\mathbb C} \backslash
G_{[a,b]}$. One has the asymptotic formula
   \[ \sqrt[s]{(a - \lambda )(b - \lambda )} \sim - \lambda \mbox { as }
      |\lambda | \rightarrow \infty .
   \]
For any $\varphi \in {\mathcal W}$ and $\lambda \in {\mathbb C} \backslash
\big( \sqcup _{k \in {\mathbb Z}} G_k \big)$ with ${\mathcal W}$ and $G_k$, $k \in {\mathbb Z}$, 
given as above, we define the {\em canonical root} of $\Delta^2(\lambda,\varphi)-4$ as
   \begin{equation}\label{10} 
        \sqrt[c]{\Delta ^2(\lambda , \varphi ) - 4} := 2i \prod_{k \in
              {\mathbb Z}}\frac{\sqrt[s]{(\lambda ^+_k(\varphi ) - \lambda )
              (\lambda ^-_k(\varphi ) - \lambda )}}{\pi _k} .
   \end{equation}
To simplify notation, we occasionally will write ${\mathcal R}(\lambda , \varphi )$
for $\Delta ^2(\lambda , \varphi ) - 4$,
   \[ {\mathcal R}(\lambda , \varphi ):= \Delta ^2(\lambda , \varphi ) - 4 .
   \]
The proof of the following lemma is straightforward and hence omitted.

\begin{Lemma}\label{Lemma 2.4} 
Let ${\mathcal W}$ be given as above. For any $\varphi \in
{\mathcal W}$, the canonical root \eqref{10} defines a holomorphic function on
${\mathbb C} \backslash (\sqcup _{k \in {\mathbb Z}} G_k)$.
\end{Lemma}
For any $\varphi \in {\mathcal W}$, define the {\em canonical sheet} (or {\em canonical branch}) 
of the open Riemann surface ${\mathcal C}^\bullet_\varphi$,
   \begin{equation}
   \label{13} {\mathcal C}^c_\varphi := \big\{ (\lambda , w) \in {\mathbb C}^2\,\big|\,
              \lambda \in {\mathbb C} \backslash (\sqcup _{k \in {\mathbb Z}}
              G_k), \ w = \sqrt[c]{\Delta ^2(\lambda , \varphi ) - 4} \}\,.
   \end{equation}
As in the Introduction, denote by $A_k$, $k\in\Z$, the cycles on the canonical sheet 
${\mathcal C}_\varphi^c$ such that for any $k\in\Z$,
\[
\pi(A_k)=\Gamma_k\,,
\] 
where $\pi : {\mathcal C}_\varphi\to\C$, $(\lambda,w)\mapsto \lambda\,$.

Finally, we need the following result on the $A_k$-periods of the holomorphic 
1-form 
$\frac{\dot\De(\lambda)}{\sqrt{\De(\lambda)^2-4}}\,d\lambda$
on ${\mathcal C}^\bullet_\varphi$.

\begin{Lemma}\label{lem:zero_form} Let $\W$ be given as above. Then, for any $\varphi\in\W$ and for any $k\in\Z$,
\[
\int_{\Ga_k}\frac{\dot\De(\lambda)}{\sqrt[c]{\De(\lambda)^2-4}}\,d\lambda\in 2\pi i \Z\,.
\]
Furthermore, for any $\psi\in L^2_c$ there exist an open neighborhood ${\cal U}$ of $\psi$ in $L^2_c$ and
an integer $N_0\ge 1$ so that for any $\varphi\in{\cal U}$ statements (i)-(iii) of 
Proposition \ref{Proposition 2.1} hold and for any $|k|\ge N_0$,
\[
\int_{\partial D_k(\pi/4)}\frac{\dot\De(\lambda)}{\sqrt[c]{\De(\lambda)^2-4}}\,d\lambda=0\,,
\]
where $\partial D_k(\pi/4)$ is the counterclockwise oriented boundary of the disk $D_k(\pi/4)$.
\end{Lemma}
\noindent{\em Proof.} Take $\varphi\in\W$ and consider the holomorphic function
$f(\lambda):=\De(\lambda)+\sqrt[c]{\De(\lambda)^2-4}$ defined for 
$\lambda\in\C\setminus\big(\sqcup_{k\in\Z}G_k\big)$. Note that $f(\lambda)$ does not vanish 
and hence the logarithm $\log f(\lambda)$ is well defined as a multi-valued function on 
$\C\setminus\big(\sqcup_{k\in\Z}G_k\big)$. Along any simple closed $C^1$-smooth curve
$\gamma : [0,1]\to\C\setminus\big(\sqcup_{k\in\Z}G_k\big)$ one can choose 
$g(t)\in\log f(\gamma(t))$ so that $g$ is continuous on $0\le t\le 1$. As $\gamma(0)=\gamma(1)$,
one has $g(1)-g(0)\in 2\pi i \Z$. In view of the identity, 
\[
d(\log f(\lambda))=\frac{\dot\De(\lambda)}{\sqrt[c]{\De(\lambda)^2-4}}\,d\lambda
\] 
the first statement of the Lemma then follows.

Let us prove the second statement of the Lemma. Take $\psi\in L^2_c$ and consider the set
of potentials $I:=\big\{t\psi\,|\,t\in[0,1]\big\}\subseteq L^2_c$.
In view of Proposition \ref{Proposition 2.1} and the compactness of the set $I$ in $L^2_c$, there
exist a connected open neighborhood ${\cal U}$ of $I$ in $L^2_c$ and an integer number $N_0\ge 1$ such that for any 
$\varphi\in{\cal U}$ and for any $|k|\ge N_0$, the statements $(i)$, $(ii)$, and $(iii)$ of 
Proposition \ref{Proposition 2.1} hold. 
Note that for any $|k|\ge N_0$, the map ${\cal U}\to\R$, $\varphi\mapsto{\mathcal P}_k(\varphi)$,
\[
{\mathcal P}_k(\varphi):=\max_{\pm}
\left\{\frac{1}{2\pi i}\int_{\partial D_k(\pi/4)}
\frac{\dot\De(\lambda,\varphi)}{\sqrt{\De(\lambda,\varphi)^2-4}}\,d\lambda\right\},
\]
where the maximum is taken over the two different choices of the square root in the integrand,
is continuous. 
As ${\mathcal P}_k(\varphi)\in\Z$ and as $\mathcal U$ is connected and $0\in{\mathcal U}$, it then
follows that ${\mathcal P}_k(\varphi)={\mathcal P}_k(0)$ for any $\varphi\in{\mathcal U}$ and for any
$|k|\ge N_0$. Finally, as  $\De(\lambda,0)=2\cos(\lambda)$ we conclude that 
$\frac{\dot\Delta(\lambda,0)}{\sqrt{\Delta^2(\lambda,0)-4}}\,d\lambda=\pm i\,d\lambda$, and hence
${\mathcal P}_k(0)=0$. This completes the proof of the Lemma.
\finishproof

\noindent Denote 
\[
\Om^*:=\frac{\dot\De(\lambda)}{\sqrt{\De(\lambda)^2-4}}\,d\lambda\,.
\]
Using the second part of Lemma \ref{lem:zero_form} we choose the neighborhood $\W$ in $L^2_\bullet$
and the integer $N_0\ge 1$ so that for any $\varphi\in\W$ and any $|k|\ge N_0$,
\begin{equation}\label{eq:zero_form1}
\int_{A_k}\Om^*=0\,.
\end{equation}
In addition, for any $|k|<N_0$,
\begin{equation}\label{eq:zero_form2}
\int_{A_k}\Om^*\in2\pi  i\,\Z,
\end{equation}
and does not depend on $\varphi\in\W$ as $\int_{A_k}\Om^*$
takes discrete values and the neighborhood $\W$ in $L^2_\bullet$ can be chosen connected.


\section{Proof of Theorems \ref{Theorem 1.0} and \ref{Theorem 1.1}}
\label{4. Proof of Theorem 1.0 and Theorem 1.1}

The aim of this section is to prove Theorem~\ref{Theorem 1.0} and
Theorem~\ref{Theorem 1.1}. Let $\W\subseteq L^2_\bullet$ be the neighborhood constructed
in Section~\ref{2. Preliminaries}. Throughout this section we fix 
$\varphi\in\W$ and define the cycles $(A_m)_{m \in {\mathbb Z}}$ as 
in Section~\ref{2. Preliminaries}. Without further reference we will use the
terminology introduced in Section~\ref{1. Introduction} and Section~\ref{2. Preliminaries}.

Let $\zeta : {\mathbb C} \rightarrow {\mathbb C}$ be entire so that $\zeta$
vanishes on the set $Z_\varphi $ of double eigenvalues of $L(\varphi )$.
It then follows that the differential 
$\omega _\zeta = \frac{\zeta (\lambda )}{\sqrt{\Delta ^2(\lambda ) - 4}}\,d\lambda $
is locally $L^2$-integrable. 
Using Stokes' theorem we will estimate 
$V(r):= \frac{i}{2} \int _{X_r} \omega _\zeta\wedge \overline{\omega _\zeta }$ 
where 
$X_r=\pi^{-1}\big(\{\lambda\in{\mathbb C}\,\big|\,|\lambda |< r\}\big) 
\cap{\mathcal C}^\bullet_\varphi$. 
In view of Proposition~\ref{Proposition 2.1} and Proposition~\ref{Proposition 2.2}
one can choose $0 < \varepsilon _m < \pi / 4, \, m > N_0$, with $\sum _{m > N_0}
\varepsilon ^2_m < \infty $ so that $\lambda ^\pm _m \in D_m(\varepsilon _m)$
and $\lambda ^\pm _{-m} \in D_{-m}(\varepsilon _m)$ for any $m > N_0$. Choose
$r > 0$ so that for some $m > N_0$
\begin{equation}\label{21} 
m\pi + \varepsilon _m \leq r \leq (m + 1) \pi - \varepsilon _{m + 1} .
\end{equation}
On the Riemann surface $X_r$ consider the points $p^+:= \pi ^{-1}(r) \cap
{\mathcal C}^c_\varphi $ and $p^-:= \pi ^{-1}(r) \backslash \{ p^+\} $ 
and choose a simple $C^1$-smooth curve $B_\ast $ on $X_r$ (i.e., a $C^1$-smooth map
$[0,1] \rightarrow X_r$ without self intersections) that connects $p^+$
with $p^-$ and changes sheets when its projection $\pi(B_*)$ passes through
the $0$-th curve $G_0$. 
If $\lambda ^+_0 = \lambda ^-_0 = \tau _0$ and hence $G_0$ is the
constant curve $\tau _0$ we allow the curve $B_*$ to pass through the point
$(\tau _0, 0) \in {\mathcal C}_\varphi $ that is excluded from $X_r$. 
The inverse image $\pi ^{-1}(\partial D_0(r))$ consists of two simple 
closed curves with images
   \begin{equation}
   \label{22} C^+_r:= \pi ^{-1} \big( \{ \lambda \in {\mathbb C} \big\arrowvert
              |\lambda | = r \} \big) \cap {\mathcal C}^c_\varphi ,
   \end{equation}
and
   \[ C^-_r:= \pi ^{-1} \big( \{ \lambda \in {\mathbb C} \big\arrowvert
              |\lambda | = r \} \big) \backslash C^+_r
   \]
that we orient so that their projections $\pi (C^\pm _r) \subseteq {\mathbb C}$
have counterclockwise orientation. Furthermore, for any $1 \leq |k| \leq m$,
denote by $B'_k$ a simple $C^1$-smooth curve in $X_r$ that starts and
ends at $p^+$ and changes sheets twice -- 
first when its projection $\pi(B_k')$ passes through $G_k$ and then through $G_0$. 
If the image of $G_k$ is a point we proceed as above and allow the curve 
$B'_k$ to pass through the point $(\tau _k, 0) \in{\mathcal C}_\varphi $. 
Similarly, for any $1 \leq |k| \leq m$, denote
by $A'_k$ a simple $C^1$-smooth curve in $X_r$ that starts and ends at
$p^+$ and that is homologous to $A_k$. The curves $B_\ast , C^+_r, A'_k$,
and $B'_k, 1 \leq |k| \leq m$, considered above, are chosen so that they
intersect each other only at $p^+$. Denote by $\tilde X_r$ the surface
obtained from $X_r$ by cutting it along the curves $B_\ast $ and 
$A'_k, B'_k$, $1 \leq |k| \leq m$. 
Then ${\tilde X}_r$ is a disk and its boundary $\partial{\tilde X}_r$
can be represented as a composition of the following curves
(composed in the order of their appearance):
$C^+_r$, $B^*$, $C^-_r$, $(B^*)^{-1}$, $B'_1$, $A'_1$, $(B'_1)^{-1}$, 
$(A'_1)^{-1}$,\,...\,,$B'_m$, $A'_m$, $(B'_m)^{-1}$, $(A'_m)^{-1}$.
Consider the function $F : {\tilde X}_r\to{\mathbb C}$ given by 
\begin{equation}\label{eq:F}
F(p):= \int ^p_{p^+} \omega _\zeta
\end{equation}
for $p \in \tilde X_r$. Note that the integral is independent of the
choice of the path and hence $F$ is well defined on $\tilde X_r$. 
Furthermore introduce $a_k = \int _{A_k} \omega _\zeta \, (0 \leq |k| \leq m)$,
$b_k := \int _{B'_k} \omega _\zeta$ $(1 \leq |k| \leq m)$, $b_\ast= \int _{B_\ast } \omega _\zeta $,
and $c^\pm _r = \int _{C^\pm _r} \omega _\zeta $. By Stokes' theorem
   \begin{align*} - 2i V(r) &= \iint _{\tilde X_r} d \big( F \overline{\omega
                     _\zeta} \big) = \int _{\partial \tilde X_r} F \overline
                     {\omega _\zeta } \\
                  &= \int _{C^+_r} F \overline {\omega _\zeta } + \int
                     _{C^-_r} F \overline {\omega _\zeta } - c^-_r \overline
                     {b_\ast } - \sum _{1 \leq |k| \leq m} \big( a_k \overline
                     {b_k} - \overline {a_k} b_k \big)
   \end{align*}
where we used that$\int _{A'_k} \omega _\zeta = a_k$ for any $1 \leq |k| \leq m$
as $A'_k$ and $A_k$ are homologous. Note that for $z^- \in C^-_r$
   \[ F(z^-) = \int _{C^+_r} \omega _\zeta + \int _{B_\ast} \omega _\zeta + \int
      ^{z^-}_{p^-} \omega _\zeta = c^+_r + b_\ast + \int ^{z^+}_{p^+} - \omega
      _\zeta
   \]
where $z^+ \in C^+_r$ is determined by $\pi (z^+) = \pi (z^-)$ and the minus
sign stems from passing to the canonical sheet. Hence
   \begin{align*} \int _{C^-_r} F(z^-) \overline {\omega _\zeta } &= - \int
                     _{C^+_r} \big( c^+_r + b_\ast - F(z^+)\big) \overline
                     {\omega _\zeta } \\
                  &= - |c^+_r|^2 - b_\ast\overline{c^+_r} + \int _{C^+_r}
                     F \overline {\omega _\zeta } 
   \end{align*}
yielding
   \[ -2 i V(r) = 
      2 \int _{C^+_r} F \overline{\omega _\zeta } -|c_r^+|^2- c^-_r\overline{b_\ast}
      - b_\ast\overline {c^+_r}
      - \sum_{1 \leq |k| \leq m} \big( a_k \overline {b_k }  -
      \overline {a_k} b_k \big) .
   \]
Now assume that $a_k = 0$ for $0 \leq |k| \leq m$. As $\sum _{0
\leq |k| \leq m} A_k$ is homologous to $C^+_r$ it then follows that $c^+
_r = 0$ and as $c^-_r = - c^+_r$ one also has $c^-_r = 0$. We thus have proved
that for any $m\pi + \varepsilon _m \leq r \leq (m + 1)\pi - \varepsilon _{m
+ 1}$ with $m > N_0$
   \begin{equation}
   \label{25} V(r) = i \int _{C^+_r} F \overline {\omega _\zeta } .
   \end{equation}
This identity will be used in the proof of the following lemma.
We recall that ${\mathcal R}(\lambda)=\Delta^2(\lambda)-4$.

\begin{Lemma}\label{Lemma 3.1} 
Let $\zeta : {\mathbb C} \rightarrow {\mathbb C}$ be entire
with $\zeta \big\arrowvert_{Z_\varphi} = 0$ and $\int _{A_j} \omega _\zeta = 0 \
\forall j \in {\mathbb Z}$. Then for any $m\pi + \varepsilon _m \leq r \leq
(m + 1) \pi - \varepsilon _{m + 1}$ with $m > N_0$
   \begin{equation}\label{26} 
         V(r) \leq \frac{1}{2} \Big( r \int ^{2\pi }_0 \Big\arrowvert
              \frac{\zeta (r e^{i\theta })}{{\mathcal R}(r e^{i\theta })}
              \Big\arrowvert d\theta \Big) ^2
   \end{equation}
and
   \begin{equation}
   \label{27} V'(r) \geq \frac{2}{r\pi} V(r) .
   \end{equation}
\end{Lemma}
\noindent{\em Proof of Lemma~\ref{Lemma 3.1}. } Using polar coordinates we get from \eqref{eq:F}
   \begin{align*} \int _{C^+_r} F \overline {\omega _\zeta } &= \int ^{2\pi }
                     _0 F(r e^{i\theta }) \overline{\Big( \frac{\zeta (r e^{i
                     \theta })}{\sqrt[c]{{\mathcal R}(r e^{i\theta })}}\,d(re^{i
                     \theta }) \Big) } \\
                   &= r^2 \int ^{2\pi }_0 \Big( \int ^\theta _0 \frac{\zeta(re
                      ^{i\theta _1})}{\sqrt[c]{{\mathcal R}(r e^{i\theta _1})}}\,
                      i e^{i\theta _1} d\theta _1 \Big) 
                      \overline{\Big(\frac{\zeta(re ^{i\theta })}
                      {\sqrt[c]{{\mathcal R}(re^{i\theta })}}
                            \Big)}
                      (-i)e^{-i\theta } d\theta
   \end{align*}
where $\sqrt[c]{{\mathcal R}(\lambda )}$ denotes the canonical root \eqref{10}.
Hence
   \begin{align*} V(r) &= \Big\arrowvert \int _{C^+_r} F(p) \overline {\omega
                     _\zeta } \Big\arrowvert\\ 
                   &\leq r^2 \int ^{2\pi }_0 \Big(
                     \int ^\theta _0 \Big\arrowvert \frac{\zeta (re^{i\theta _1})}
                     {\sqrt[c]{{\mathcal R}(re^{i\theta _1})}} \Big\arrowvert
                     d\theta _1 \Big) d \Big(
                     \int ^\theta _0 \Big\arrowvert \frac{\zeta (re^{i\theta_1})}
                     {\sqrt[c]{{\mathcal R}(re^{i\theta_1})}} \Big\arrowvert
                     d\theta_1 \Big) \\
                  &= \frac{r^2}{2} \Big( \int ^{2\pi }_0 \Big\arrowvert
                     \frac{\zeta (re^{i\theta _1})}
                     {\sqrt[c]{{\mathcal R}(re^{i\theta _1})}} \Big\arrowvert
                     d\theta _1 \Big) ^2 .
   \end{align*}
This proves the estimate \eqref{26}. To get \eqref{27} note that
   \begin{align*} &\frac{r^2}{2} \Big( \int ^{2\pi }_0 \Big\arrowvert
                     \frac{\zeta (re^{i\theta })}
                     {\sqrt[c]{{\mathcal R}(re^{i\theta })}} \Big\arrowvert
                     d\theta \Big) ^2 \leq \pi r^2 \int ^{2\pi } _0
                     \Big\arrowvert \frac{\zeta (re^{i\theta })}
                     {\sqrt[c]{{\mathcal R}(re^{i\theta })}} \Big\arrowvert ^2
                     d\theta  \\
                  &= \pi r \frac{d}{dr} \Big( \int ^r_0 \int ^{2\pi }_0 \Big
                     \arrowvert \frac{\zeta (\rho e^{i\theta })}
                     {\sqrt[c]{{\mathcal R}(\rho e^{i\theta })}} \Big
                     \arrowvert ^2 \rho\,d\theta d\rho\Big) 
                     = \pi r \frac{V'(r)}{2} .
   \end{align*}
Hence $V(r) \leq \frac{r\pi}{2} V'(r)$ as claimed.
\finishproof

Estimate \eqref{27} is now used to prove Theorem~\ref{Theorem 1.0}.

\medskip

\noindent{\em Proof of Theorem~\ref{Theorem 1.0}. } Assume that $\zeta \not\equiv
0$ is an entire function satisfying $\int _{A_m} \omega _\zeta = 0$ for any $m
\in {\mathbb Z}$. Then $V(r) \not= 0 \ \forall r > 0$ and in view of \eqref{27},
for any $m>N_0$ and for any $m\pi+\varepsilon_m\leq r\leq(m+1)\pi-\varepsilon_{m+1}$,
\[ 
(\log V(r))' \geq \frac{2}{\pi r}\,.
\]
Integrating this inequality over the interval $[m\pi+\varepsilon_m,(m+1)\pi-\varepsilon_{m+1}]$
we obtain that
   \[ \frac{V\big((m + 1)\pi - \varepsilon _{m+1}\big)}{V(m\pi + \varepsilon _m)}
      \geq e^{\frac{2}{\pi }\big( \log ((m + 1)\pi - \varepsilon _{m+1}) -
      \log(m\pi + \varepsilon _m) \big)} .
   \]
This implies that for any $m \geq m_0 > N_0$,
   \[ V \big( (m + 1) \pi \big) \geq V(m_0 \pi ) e^{\frac{2}{\pi } S(m,m_0)},
   \]
where
   \begin{align*} &S(m,m_0) := \sum ^m_{j=m_0} \log \big( (j+1) \pi -
                     \varepsilon _{j + 1}\big) - \log (j \pi + \varepsilon _j) \\
                  &= \log \frac{m+1}{m_0} + \sum ^m_{j = m_0} \log \Big( 1 -
                     \frac{\varepsilon _{j + 1}}{(j + 1)\pi }\Big) - \sum ^m_{j =
                     m_0} \log \Big( 1 + \frac{\varepsilon _j}{j\pi } \Big) \\
                  &\geq \log \frac{m + 1}{m_0} - O \Big( \sum ^{m+1}_{j=m_0}
                     \frac{\varepsilon _j}{j\pi } \Big) .
   \end{align*}
As $\sum ^m _{j=m_0} \frac{\varepsilon _j}{j\pi } = O \Big( \sum ^m_{j=m_0}
\varepsilon ^2_j \Big) ^{1/2}$ one then concludes that
   \[ V\big((m + 1)\pi\big) \geq C( m + 1)^{2/\pi }
   \]
where $C > 0$ depends on $m_0 > N_0$ but not on $m \geq m_0$. 
\finishproof

\noindent{\em Proof of Theorem~\ref{Theorem 1.1}. } Theorem~\ref{Theorem 1.1}
is an immediate consequence of Theorem~\ref{Theorem 1.0}.
\finishproof


\section{Outline of proof of Theorem~\ref{Theorem 1.2}}
\label{3. Outline of proof of Theorem 1.2}

In this section we describe the main steps in the proof of Theorem~\ref{Theorem 1.2}. 
For the remaining part of the paper $\W$ will denote the neighborhood
${\mathcal W} \subseteq L^2_\bullet$ constructed in Section~\ref{2. Preliminaries} with $N_0 \geq 1$ so that
Proposition \ref{Proposition 2.1}, Proposition \ref{Proposition 2.2}, and the identities
\eqref{eq:zero_form1} and \eqref{eq:zero_form2} hold.
Let $(A_m)_{m \in {\mathbb Z}}$ be the cycles on the canonical branch of 
${\mathcal C}_\varphi^\bullet$ introduced in Section~\ref{2. Preliminaries}.
Without further explanations, for any given $\varphi \in {\mathcal W}$ and
$n \in {\mathbb Z}$ consider the following ansatz for the holomorphic differentials
of Theorem~\ref{Theorem 1.2}
   \begin{equation}\label{15} 
     \Omega ^n_\beta := \Omega ^n - \omega ^n_\beta
   \end{equation}
where the forms $\Omega ^n$ and $\omega ^n_\beta$ are defined as follows:
   \begin{equation}\label{16} 
           \Omega ^n := \begin{cases} \frac{\dot \Delta (\lambda )}{\lambda
              - \dot \lambda _n} \frac{d\lambda }{\sqrt{{\mathcal R}(\lambda )}} ,
              &\,|n| > N_0 \\ \frac{\dot \Delta (\lambda )}{\lambda
              - \dot \lambda _{N_0}} \frac{d\lambda }{\sqrt{{\mathcal R}(\lambda )}} ,
              &\,|n| \leq N_0 \end{cases}
   \end{equation}
and for any given $\beta = (\beta _k)_{k \not= n} \in \ell ^1_{\hat n} \equiv
\ell ^1_{\hat n,{\mathbb C}} := \ell ^1\big({\mathbb Z} \backslash \{ n \},{\mathbb C}\big)$
\begin{equation}\label{17}
\omega^n_\beta := \frac{\xi ^n_\beta (\lambda )}{\sqrt{{\mathcal R}(\lambda)}}\,d\lambda
\end{equation}
where
\begin{align}\label{17uno} 
                 &\xi ^n_\beta (\lambda ):= \begin{cases} \Big( \sum\limits_
                   {|j| > N_0,j\ne n} \frac{\beta _j}{\lambda - \dot \lambda
                   _j} + \frac {p^n_\beta (\lambda )}{\prod _{|k| \leq N_0}
                   (\lambda - \dot \lambda _k)} \Big) \frac{\dot \Delta (\lambda )}
                   {\lambda - \dot \lambda _n} , &\, |n| > N_0 \\
                   \Big( \sum\limits_{|j| > N_0} \frac{\beta _j}{\lambda - \dot \lambda
                   _j} + \frac{p^n_\beta (\lambda )}{\prod _{\underset{k \not=
                   N_0}{|k| \leq N_0}} (\lambda - \dot \lambda _k)}\Big)
                   \frac{\dot \Delta (\lambda )}{\lambda - \dot \lambda _{N_0}} ,
                   &\,|n| \leq N_0 \end{cases}
  \end{align}
and
   \[ p^n_\beta (\lambda ) := \begin{cases} \,\,\,\,\,\sum\limits^{2N_0}_{j=0} \beta _{j - N_0}
      \lambda ^j , &\, |n| > N_0 \\ \sum\limits_{j=0}^{n+N_0-1} \beta
      _{j - N_0} \lambda ^j + \sum\limits_{j = n + N_0}^{2N_0 - 1} \beta _{j - N_0 + 1}
      \lambda ^j , &\,|n| \leq N_0 . \end{cases}
   \]
Note that $\xi ^n_\beta (\lambda )$ is entire. In the case $|n| \leq N_0$, it
is convenient to write the polynomial $p^n_\beta (\lambda )$ in the following
alternative way 
$p^n_\beta (\lambda ) = \big( \sum_{|j|\le N_0, j\ne n} 
\beta _j \lambda ^{j - \varepsilon ^n_j}\big) \lambda ^{N_0}$ 
where
   \begin{equation}\label{eq:epsilon}
      \varepsilon ^n_j := \begin{cases} 1 , &j \geq n \\ 0 ,
      &\, j < n .\end{cases}
   \end{equation}
We want to find $\beta = (\beta _k)_{k \not= n} \in \ell ^1_{\hat n}$ so that
   \begin{equation}
   \label{17bis} \frac{1}{2\pi } \int _{A_m} \Omega ^n_\beta = \delta _{nm}
                 \quad \forall m \in {\mathbb Z} .
   \end{equation}

\noindent The following proposition is proved in Section~\ref{6. Proof of Prop 4.1}.

\begin{Prop}\label{Proposition 4.1} 
For any $n \in {\mathbb Z}$, $\varphi \in {\mathcal W}$, and $\beta \in \ell ^1_{\hat n}$
   \begin{equation}
   \label{17ter} \sum _{m \in {\mathbb Z}} \frac{1}{2\pi } \int _{A_m} \Omega
                 ^n_\beta = \lim _{K \rightarrow \infty } \sum _{|m| \leq K}
                 \frac{1}{2\pi } \int _{A_m} \Omega ^n_\beta = 1 .
   \end{equation}
In particular, for $\beta = 0$, $\Omega^n_\beta\equiv\Omega^n$ satisfies
   \[ \sum _{m \in {\mathbb Z}} \frac{1}{2\pi } \int _{A_m} \Omega ^n =
       1\,.
   \]
\end{Prop}
In view of Proposition~\ref{Proposition 4.1}, the system of equations
\eqref{17bis} is equivalent to
   \begin{equation}
   \label{17quater} \int _{A_m} \omega ^n_\beta = \int _{A_m} \Omega ^n \quad\quad
                    \forall m \not= n .
   \end{equation}
By multiplying the right and left hand side of the above equation by 
$m\pi-n\pi$ (if $|n| > N_0$) or $m\pi-N_0\pi $ (if $|n| \leq N_0$), we arrive
at the following linear system for $\beta $,
   \begin{equation}
   \label{18} T^n \beta = b^n,\,\,\,\,\, T^n = (T^n_{mj})_{m,j \not= n},\,\,\,\,
                           b^n = (b^n_m)_{m \not= n}
   \end{equation}
where for any $m \not= n$, $b^n_m$ is given by
   \[ b^n_m := \begin{cases} \frac{1}{2\pi } \int _{\Gamma _m} \frac{m\pi -
      n\pi }{\lambda - \dot \lambda _n} \frac{\dot \Delta (\lambda )}{\sqrt[c]
      {{\mathcal R}(\lambda )}}\,d \lambda , &\, |n| > N_0 \\ 
      \frac{1}{2\pi } \int _{\Gamma _m} \frac{m\pi - N_0 \pi }{\lambda - \dot \lambda
      _{N_0}} \frac{\dot \Delta (\lambda )}{\sqrt[c]{{\mathcal R}(\lambda )}}\,
      d\lambda , &\, |n| \leq N_0 \end{cases}
   \]
and for $|n| > N_0$, $T^n_{mj}$ is given by
   \[ T^n_{mj} := \begin{cases} \frac{1}{2\pi } \int _{\Gamma _m} \frac{m\pi -
      n\pi }{\lambda - \dot \lambda _n} \frac{1}{\lambda - \dot \lambda _j}
      \frac{\dot \Delta (\lambda )}{\sqrt[c]
      {{\mathcal R}(\lambda )}}\,d \lambda , &\, |j| > N_0,\, j\ne n \\ 
      \frac{1}{2\pi } \int _{\Gamma _m} \frac{m\pi - n \pi }{\lambda - \dot \lambda
      _{n}} \frac{\dot \Delta (\lambda )}{\prod _{|k| \leq N_0}(\lambda
      - \dot \lambda _k)}\frac{\lambda ^{N_0 + j}}{\sqrt[c]{{\mathcal R}(\lambda )}}\,
      d\lambda , &\, |j| \leq N_0 \end{cases}
   \]
whereas for $|n| \leq N_0$ one has
   \[ T^n_{mj} := \begin{cases} \frac{1}{2\pi } \int _{\Gamma _m} \frac{m\pi -
      N_0\pi }{\lambda - \dot \lambda _{N_0}} \frac{1}{\lambda - \dot \lambda _j}
      \frac{\dot \Delta (\lambda )}{\sqrt[c]
      {{\mathcal R}(\lambda )}}\,d \lambda , &\, |j| > N_0 \\ 
      \frac{1}{2\pi}\int_{\Gamma_m}
      \frac{(m\pi-N_0\pi)\dot\Delta(\lambda)}{\prod_{|k|\leq N_0}(\lambda-\dot\lambda_k)}
      \frac{\lambda^{N_0+j-\varepsilon^n_j}}{\sqrt[c]{{\mathcal R}(\lambda )}}\,
      d\lambda, &\, |j| \leq N_0,\, j\ne n . \end{cases}
   \]
Using Proposition~\ref{Proposition 6.2} -- an application of
Theorem~\ref{Theorem 1.1} -- we prove in Section~\ref{7. Existence of normalized
differentials} the following

\begin{Prop}\label{Proposition 4.2} 
For any $n \in {\mathbb Z}$ and $\varphi \in{\mathcal W}$ we have: 
\begin{itemize}
\item[(i)] $b^n \in \ell ^1_{\hat n}$;
\item[(ii)] $T^n : \ell ^1_{\hat n}\rightarrow \ell ^1_{\hat n}$ is a linear isomorphism.
\end{itemize}
\end{Prop}
Denote by $\beta^n\equiv\beta^n(\varphi)\in\ell^1_{\hat n}$ the unique
solution of \eqref{18}, guaranteed by Proposition~\ref{Proposition 4.2} and define,
\begin{equation}\label{eq:zeta_n}
\zeta _n(\lambda , \varphi ) := \begin{cases} \frac{\dot \Delta (\lambda )}
      {\lambda - \dot \lambda _n} - \xi ^n_{\beta ^n}(\lambda ) , &
      |n| > N_0 \\ \frac{\dot \Delta (\lambda )}{\lambda - \dot \lambda _{N_0}}
      - \xi ^n_{\beta ^n} (\lambda ) , & |n| \leq N_0, \end{cases}
\end{equation}
where $\xi ^n_{\beta ^n}$ is given by \eqref{17uno} with $\beta ^n$ substituted
for $\beta $. The following proposition is proved in 
Section~\ref{7. Existence of normalized differentials}.
\begin{Prop}\label{Proposition 4.3} 
For any $n \in {\mathbb Z}$, $\beta ^n : {\mathcal W} \rightarrow \ell ^1_{\hat n}$ and
$\zeta _n : {\mathbb C}\times {\mathcal W} \rightarrow {\mathbb C}$ are analytic maps.
Furthermore, for any $\varphi \in {\mathcal W}$ and $n \in {\mathbb Z}$,
   \[ \frac{1}{2\pi } \int _{A_m} \frac{\zeta _n(\lambda , \varphi )}{\sqrt{
      {\mathcal R}(\lambda , \varphi )}}\,d \lambda = \delta _{nm} \quad \forall
      m \in {\mathbb Z} .
   \]
\end{Prop}

To obtain uniform in $n\in\Z$ and locally uniform in $\W$ estimates of the zeros of $\zeta _n$ we
consider the following ``limiting'' linear system for 
$\beta=(\beta_k)_{k \in {\mathbb Z}} \in \ell ^1 \equiv \ell ^1_{\mathbb C}$,
   \begin{equation}\label{19}
     T^\ast \beta = b^\ast
   \end{equation}
where $T^\ast = (T^\ast _{mj})_{m,j \in {\mathbb Z}}$ is given by
   \[ T^\ast _{mj} := \begin{cases} \frac{1}{2\pi } \int _{\Gamma _m} \frac{1}
      {\lambda - \dot \lambda _j} \frac{\dot \Delta (\lambda )}{\sqrt[c]{
      {\mathcal R}(\lambda )}}\,d\lambda , &\, |j| > N_0 \\ \frac{1}
      {2\pi } \int _{\Gamma _m} \frac{\lambda ^{N_0 + j}}{\prod _{|k| \leq N_0}
      (\lambda - \dot \lambda _k)} \frac{\dot \Delta (\lambda )}{\sqrt[c]
      {{\mathcal R}(\lambda )}}\, d \lambda , &\, |j| \leq N_0
      \end{cases}
   \]
and $b^\ast:=(b^\ast_m)_{m\in\Z}$ is given by
\[
b^\ast_m:=\frac{1}{2\pi } \int _{\Gamma _m}
           \frac{\dot \Delta (\lambda )}{\sqrt[c]{{\mathcal R}(\lambda )}}\,d\lambda\,. 
\]
This linear system is equivalent to the condition
\begin{equation}\label{eq:linear_system*}
\int_{A_m}\Om^*_\beta=0\,\,\,\,\forall m\in\Z
\end{equation}
where $\beta \in \ell ^1$ and the holomorphic $1$-form $\Om^*_\beta$ on
$C^\bullet_\varphi $ is given by
\[
\Om^*_\beta:=\Om^*-\om^*_\beta
\]
where
\[
\om_\beta^*:=\frac{\xi ^\ast _\beta(\la)}{\sqrt{{\mathcal R}(\la )}}\,d\la
\]
with
\begin{equation}\label{eq:zeta_beta}
\xi^\ast _\beta(\la):=\Big(\sum_{|j|>N_0}\frac{\beta_j}{\lambda - 
\dot \lambda _j}+\frac{p^\ast _\beta(\la)}{\prod_{|j|\le N_0}(\la-\dot \lambda_j)}\Big)
\dot \Delta (\la)
\end{equation}
and
\[
p^\ast_\beta(\la):= \sum_{j=0}^{2N_0} \beta _{-N_0 + j} \lambda ^j .
\]
Note that $\xi ^\ast _\beta $ is an entire function of $\lambda$ and
$p^\ast_\beta(\lambda)$ is a polynomial of degree at most $2N_0$.
We can rewrite \eqref{eq:linear_system*} as
\[
\int_{A_m}\om^*_\beta=\int_{A_m}\Om^*\,\,\,\forall m\in\Z
\]
that leads to the linear system \eqref{19} in view of the definitions of
the forms $\om^*_\beta$ and $\Om^*$.
Using Proposition~\ref{Proposition 5.1} -- another application of
Theorem~\ref{Theorem 1.1} -- we prove in 
Section~\ref{7. Existence of normalized differentials} the following
\begin{Prop}\label{Proposition 4.4}
For any $\varphi \in {\mathcal W}$, $
T^\ast : \ell ^1\rightarrow \ell ^1$ is a linear isomorphism.
\end{Prop}
Recall that in view of \eqref{eq:zero_form1} and \eqref{eq:zero_form2},
$b^*\equiv b^*(\varphi)$ is in $\ell^1$ and does not depend on $\varphi\in\W$.
Denote by $\beta^*\equiv\beta^*(\varphi)\in\ell^1$ the unique solution of 
\eqref{19} guaranteed by Proposition \ref{Proposition 4.4}.

Proposition~\ref{Proposition 4.4} is used in 
Section~\ref{Estimates of the zeros} to prove the uniform estimates of 
the zeros of $\zeta _n$, stated in Theorem~\ref{Theorem 1.2} -- 
see Proposition~\ref{Proposition 8.3} and Lemma~\ref{Lemma 8.5} -- 
and the product representation of $\zeta _n $ --
see Corollary~\ref{Corollary 8.4}. Finally, in Lemma~\ref{Lemma 8.6}, it is
proved that $\zeta _n(\cdot , \varphi )$ vanishes on the set 
$Z_\varphi \backslash\{ \lambda ^\pm _n (\varphi )\}$. Combining the results described above,
the proof of Theorem~\ref{Theorem 1.2} is complete.


\section{Vanishing Lemma}\label{sec:application_of_vl}

Let us fix $\varphi\in\W\subseteq L^2_\bullet$ where as stated in the beginning of 
Section \ref{3. Outline of proof of Theorem 1.2} $\W$ denotes the neighborhood constructed in 
Section~\ref{2. Preliminaries}. In this section we prove the following

\begin{Prop}\label{Proposition 5.1} 
Let $\beta \in \ell ^1$ be arbitrary. If $\int _{A_m}\omega _\beta^* = 0$
for any $m \in {\mathbb Z}$, then $\beta = 0$.
\end{Prop}
We prove Proposition~\ref{Proposition 5.1} with the help of Theorem~\ref{Theorem 1.1}.
To this end we prove the following lemmas.
\begin{Lemma}\label{Lemma 5.2}
If $\int_{A_m}\om_\beta^*=0$ for any $m\in\Z$ then $\xi ^\ast _\beta \big\arrowvert_{Z_\varphi} = 0$.
\end{Lemma}
\noindent{\em Proof of Lemma~\ref{Lemma 5.2}.} Assume that for some $k\in{\mathbb Z}$, 
$\lambda ^-_k = \lambda ^+_k = \tau _k$. Then in view of \eqref{10} for 
$(\lambda,w)\in{\mathcal C}^c_\varphi$ near $(\tau_k,0)\in{\mathcal C}_\varphi$,
   \begin{equation}
   \label{29} \omega ^\ast _\beta = \xi ^\ast _\beta (\lambda ) \frac{h(\lambda )}
              {\lambda - \tau _k}\, d\lambda
   \end{equation}
where $h(\lambda )$ is a holomorphic function that is defined in an open
neighborhood of $\tau _k$ and satisfies $h(\tau _k) \not= 0$. As by
assumption $\int_{A_k} \omega ^\ast _\beta = 0$ we conclude from
\eqref{29} that $\xi ^\ast _\beta (\tau _k) = 0$.
\finishproof

Lemma~\ref{Lemma 5.2} implies that for any $\beta \in \ell ^1$ as in
Proposition~\ref{Proposition 5.1} and for any $r > 0$
   \[ V_\beta (r):= \frac{i}{2} \int _{X_r} \omega ^\ast _\beta \wedge
      \overline{\omega ^\ast _\beta } < \infty .
   \]

\begin{Lemma}\label{Lemma 5.3} 
If $\int _{A_j} \omega ^\ast _\beta = 0$ for any $j \in{\mathbb Z}$, then for any $\delta > 0$,
   \begin{equation}\label{30} 
              V_\beta (m\pi ) = O(m^\delta ) \ \mbox { as } m \rightarrow\infty .
   \end{equation}
\end{Lemma}
\noindent{\em Proof of Lemma~\ref{Lemma 5.3}. } Assume that $\int _{A_j}\omega^\ast_\beta = 0$ for any 
$j \in {\mathbb Z}$. Then by \eqref{26} in Lemma~\ref{Lemma 3.1}, 
for $r_m = (m + \frac{1}{2}) \pi $ with $m > N_0$,
   \[ V_\beta (m\pi ) \leq \frac{r^2_m}{2} \left( \int ^{2\pi }_0 \left|
      \frac{\xi ^\ast _\beta (r_m e^{i\theta })}{\sqrt[c]{{\mathcal R}(r_m
      e^{i\theta })}} \right| d\theta \right) ^2 .
   \]
By \cite[Lemma C.5]{GKP} and Proposition~\ref{Proposition 2.2},
   \[ \frac{\dot \Delta (r_m e^{i\theta })}{\sqrt[c]{{\mathcal R}(r_m e^{i
      \theta })}} = O (1) \ \mbox{ as } m \rightarrow \infty
   \]
uniformly in $0 \leq \theta < 2\pi $. This together with the definition of
$\xi ^\ast _\beta $ implies
   \begin{equation}
   \label{31} \left| \frac{\xi ^\ast _\beta (r_m e^{i \theta }) }{\sqrt
              [c]{{\mathcal R}(r_m e^{i\theta })}}\right| \leq C
              \left( \sum _{|j| > N_0} \frac{|\beta _j|}{|r_m e^{i\theta } -
              \dot \lambda _j|} + \frac{|p^\ast _\beta (r_m e^{i\theta })|}
              {\prod _{|j| \leq N_0} |r_m e^{i\theta } - \dot \lambda _j|}\right)
   \end{equation}
with a constant $C > 0$ independent of $m > N_0$ and $0 \leq \theta < 2\pi $.
As $p^\ast _\beta (\lambda )$ is a polynomial in $\lambda $ of degree at most
$2N_0$ it follows that
   \begin{equation}
   \label{31bis} \int ^{2\pi }_0 \frac{|p^\ast _\beta (r_m e^{i\theta })|}
                 {\prod _{|j| \leq N_0} |r_m e^{i\theta } - \dot \lambda
                 _j|}\,d\theta = O \Big( \frac{1}{m} \Big) \ \mbox{ as }
                 m \rightarrow \infty .
   \end{equation}
For any $m > 2N_0$, we split the sum $\sum _{|j| > N_0} = \sum _{j \in J_1(m)}
+ \sum _{j \in J_2(m)}$ where $J_1(m)$ is the set 
   \[ \Big\{ j > N_0 \Big\arrowvert \Big| j - \Big( m + \frac{1}{2}\Big) \Big| 
      \leq \frac{m}{2} \Big\} \cup \Big\{ j < - N_0 \Big\arrowvert
      \Big| j + \Big( m + \frac{1}{2} \Big) \Big| \leq \frac{m}{2} \Big\}
   \]
and $J_2(m):= {\mathbb Z} \backslash ([ - N_0, N_0] \cup J_1(m))$. Then for
any $j \in J_2(m)$ with $m > 2N_0$ and any $0 \leq \theta < 2\pi $
   \[ |r_m e^{i\theta } - \dot \lambda _j| \geq \Big\arrowvert \Big( m +
      \frac{1}{2} \Big) \pi - j\pi \Big\arrowvert - \frac{\pi }{4} \geq
      \frac{(m+1)\pi }{2}-\frac{\pi}{4}\ge\frac{m\pi}{2}
   \]
implying that
   \begin{equation}
   \label{32} \int ^{2\pi }_0 \sum _{j \in J_2(m)} \frac{|\beta _j|}{|r_m e
              ^{i\theta   } - \dot \lambda _j|}\,d\theta = O \Big( \frac{1}{m} \Big) .
   \end{equation}
To estimate $\int ^{2\pi }_0 \sum _{j \in J_1(m)} \frac{|\beta _j|}{|r_m e^{i
\theta } - \dot \lambda _j|}\,d\theta $ the integral $\int ^{2\pi }_0$ is
split up as follows: For $0 < \alpha < 1$, one has uniformly in $j \in J_1(m)$,
   \[ \int ^{\frac{1}{m^\alpha }} _{-\frac{1}{m^\alpha }} \frac{d\theta }
      {|r_m e^{i \theta } - \dot \lambda _j|} , \ \int ^{\pi + \frac{1}{m^\alpha }}
      _{\pi -\frac{1}{m^\alpha }} \frac{d\theta }{|r_m e^{i \theta } - \dot \lambda
      _j|} = O \Big( \frac{1}{m^\alpha } \Big)
   \]
as $|r_m e^{i\theta}-\dot\lambda_j|\ge\pi/4$.
By choosing $N_1 \geq 2N_0$ sufficiently large we can ensure that for any $m \geq
N_1$, $\theta \in \big[ \frac{1}{m^\alpha } , \pi - \frac{1}{m^\alpha }\big]$ and
$j \in J_1(m)$
   \begin{align}
   \begin{split}
   \label{33} |r_m e^{i\theta } - \dot \lambda _j| &\geq |r_m \sin \theta -
                \mathop{\rm Im}{\dot\lambda _j}| \\
              &\geq r_m \sin \theta - \frac{\pi }{4} \geq \Big( \Big( m +
                 \frac{1}{2} \Big) \sin \frac{1}{m^\alpha } - \frac{1}{4}
                 \Big) \geq 1 .
   \end{split}
   \end{align}
Note that $r_m \sin \theta - \frac{\pi }{4}$ is the distance from $r_m e^{i\theta }$ to 
the horizontal line $\mathop{\rm Im} z = \frac{\pi }{4}$. Using \eqref{33},
$\sin \theta \geq \frac{2}{\pi } \theta $ for $0 \leq \theta \leq \pi / 2$
and taking $N_1 \geq 2N_0$ larger if necessary we get
   \begin{align*} &\int ^{\pi - \frac{1}{m^\alpha }} _{\frac{1}{m^\alpha }}
                     \frac{d\theta }{|r_m e^{i \theta } - \dot \lambda _j|} 
                     \leq 2 \int ^{\frac{\pi }{2}}_{\frac{1}{m^\alpha }}
                     \frac{d\theta }{r_m \sin \theta   - \frac{\pi }{4}} \\
                  &= \frac{\pi}{r_m} \int ^{\pi / 2}_{\frac{1}{m^\alpha }}
                     \frac{d\theta }{\theta - \pi ^2 / 8 r_m} = O\Big(\frac{\log m}{m}\Big)
   \end{align*}
and similarly
   \[ \int ^{2\pi - \frac{1}{m^\alpha }} _{\pi +\frac{1}{m^\alpha }}
      \frac{d\theta }{|r_m e^{i \theta } - \dot \lambda _j|} = O\Big(\frac{\log m}{m}\Big)
   \]
uniformly in $j \in J_1(m)$. Hence
   \begin{equation}
   \label{34} \int ^{2\pi }_0 \sum _{j \in J_1(m)} \frac{|\beta _j|}{|r_m e^{i
              \theta } - \dot \lambda _j|} d\theta = O\Big(\frac{1}{m^\alpha}\Big) \ 
              \mbox { as } m \rightarrow \infty .
   \end{equation}
Combining \eqref{31}, \eqref{31bis}, \eqref{32} and \eqref{34} yields
   \[ V_\beta (m\pi ) = O\big( m^{2 - 2\alpha } \big)
   \]
for any $0 < \alpha < 1$. This completes the proof of the lemma.
\finishproof

\noindent{\em Proof of Proposition~\ref{Proposition 5.1}. } By assumption
$\omega ^\ast _\beta $ with $\beta \in \ell ^1_{\mathbb C}$ satisfies
$\int _{A_m} \omega ^\ast _\beta = 0$ for any $m \in {\mathbb Z}$. Then
Lemma~\ref{Lemma 5.2} implies that $\xi ^\ast _\beta \big\arrowvert _{Z_\varphi }
= 0$. Hence we can apply Theorem~\ref{Theorem 1.1} and Lemma~\ref{Lemma 5.3} to
conclude that $\xi ^\ast _\beta \equiv 0$. Evaluating $\xi ^\ast _\beta $
at $\lambda = \dot \lambda _k$ with $|k| > N_0$ we get from \eqref{eq:zeta_beta} and
Proposition \ref{Proposition 2.3} that
   \[ 0 = \xi ^\ast _\beta (\dot \lambda _k) = 2 \beta _k \prod _{m \not=
      k} \frac{\dot \lambda _k - \dot \lambda _m}{\pi _m} .
   \]
As $|k| > N_0$, $\dot\lambda_k$ is a simple zero of $\dot \Delta (\lambda )$
(cf. Proposition~\ref{Proposition 2.1}) and hence $\prod _{m \not= k}
\frac{\dot \lambda _k - \dot \lambda _m}{\pi _m} \not= 0$. We therefore
conclude that $\beta _k = 0$ for any $|k| > N_0$ and thus in view of \eqref{eq:zeta_beta},
   \[ \xi ^\ast _\beta (\lambda ) = p^\ast _\beta (\lambda )\,\frac{\prod
      _{|m| > N_0} \frac{\lambda - \dot \lambda _m}{\pi _m}}{\prod _{|m|
      \leq N_0} \pi _m } .
   \]
As $\xi ^\ast _\beta \equiv 0$ it follows that $p^\ast _\beta \equiv 0$
implying that $\beta _k = 0$ for $|k| \leq N_0$. We thus have proved that
$\beta = 0$ as claimed.
\finishproof

For any given $n \in {\mathbb Z}$ and $\beta = (\beta _k)_{k\not= n} \in
\ell ^1_{\hat n}$ consider the holomorphic $1$-form
$\frac{\xi ^n_\beta(\lambda )}{\sqrt{{\mathcal R}(\lambda )}}\,d\lambda $ defined
in \eqref{17}. Arguing in the same way as in the proof of 
Proposition~\ref{Proposition 5.1} one obtains

\begin{Prop}
\label{Proposition 5.4} Let $n \in {\mathbb Z}$ and $\beta = (\beta _k)_{k
\not= n} \in \ell ^1_{\hat n}$ be arbitrary. If $\int _{A_m} \omega ^n_\beta
= 0$ for any $m \in {\mathbb Z}$, then $\beta = 0$.
\end{Prop}


\section{Proof of Proposition~\ref{Proposition 4.1}}\label{6. Proof of Prop 4.1}

The aim of this section is to prove Proposition~\ref{Proposition 4.1} concerning
the identity of the sum of all $A$-periods of a holomorphic differential of the
form $\Omega ^n_\beta $.

\noindent{\em Proof of Proposition~\ref{Proposition 4.1}. } As the proof in
the two cases $|n| \leq N_0$ and $|n| > N_0$ are similar, we consider the case
$|n| > N_0$ only. Recall that for any $|n| > N_0$ and $\beta\in\ell ^1_{\hat n}$
   \[ \Omega ^n_\beta = \frac{1 - \eta ^n_\beta (\lambda )}{\lambda - \dot
      \lambda _n} \frac{\dot \Delta (\lambda )}{\sqrt{\mathcal R(\lambda )}}\,d \lambda
   \]
where
   \begin{equation}
   \label{35} \eta ^n_\beta (\lambda ) = \sum _{|j| > N_0,j\ne n}
              \frac{\beta _j}{\lambda - \dot \lambda _j} + \frac{p^n_\beta
              (\lambda )}{\prod _{|k| \leq N_0} (\lambda - \dot \lambda _k)}
   \end{equation}
is a meromorphic function that can have poles only at the points $\dot
\lambda _j, j \in {\mathbb Z}$. Let $r_m:= m\pi + \pi / 2$ for $m > N_0$
with $N_0$ as in Proposition~\ref{Proposition 2.1}. As $\sum _{|k| \leq
m} A_k$ is homologous to $C^+_r$ (cf. \eqref{22}) one has for any $m > N_0$
   \begin{equation}
   \label{36} \sum _{|k| \leq m} \frac{1}{2\pi } \int _{A_k} \Omega ^n
              _\beta = \frac{1}{2\pi } \int _{C^+_{r_m}} \Omega ^n_\beta .
   \end{equation}
In view of \cite[Lemma C.5]{GKP}, Proposition~\ref{Proposition 2.1},
Proposition~\ref{Proposition 2.2}, and the definition of the canonical
root \eqref{10}
   \[ \frac{\dot \Delta (r_m e^{i\theta })}{\sqrt[c]{{\mathcal R}(r_m e
      ^{i\theta })}} = \frac{1}{i} (1 + o(1))\,\,\,\mbox { as } m \rightarrow
      \infty
   \]
uniformly for $0 \leq \theta < 2\pi $. This together with \eqref{35} implies
that uniformly for $\lambda \in C^+_{r_m}$
   \begin{equation}
   \label{37} \Omega ^n_\beta = \frac{1}{i} \Big( \frac{1}{\lambda - \dot
              \lambda _n} + O \Big( \frac{\eta ^\beta _n(\lambda )}{m} \Big)
              + o \Big( \frac{1}{m} \Big) \Big)\,d\lambda \,\mbox { as } m
              \rightarrow \infty
   \end{equation}
with constants independent of $\lambda \in C^+_{r_m}$ and $m > N_0$.
Combining \eqref{36} with \eqref{37} one gets for 
$m>\max\{n,N_0\} $
   \begin{equation}
   \label{38} \sum _{|k| \leq m} \frac{1}{2\pi } \int _{A_k} \Omega ^n_\beta
              = 1 + O \big( \max _{C^+_{r_m}} |\eta ^\beta _n (\lambda )| \big)
              + o(1)
   \end{equation}
with constants uniform in $m>\max\{n,N_0\}$. Arguing in a similar way
as in the proof of Lemma~\ref{Lemma 5.3} one sees that
   \[ \max _{\lambda \in C^+_{r_m}} |\eta ^\beta _n(\lambda ) | \rightarrow 0
      \mbox { as } m \rightarrow \infty .
   \]
This combined with \eqref{38} yields Proposition~\ref{Proposition 4.1}.
\finishproof

As an immediate Corollary of Proposition~\ref{Proposition 4.1} we get the
following result for 
$\omega^n_\beta=\Omega^n-\Omega^n_\beta={\Omega^n_\gamma}\big|_{\gamma=0}-\Omega^n_\beta$.

\begin{Coro}\label{Corollary 6.1} 
For any $n \in {\mathbb Z}$ and any $\beta \in \ell ^1
_{\hat n}$, 
   \[ \sum _{m \in {\mathbb Z}} \frac{1}{2\pi }
      \int _{A_m} \omega ^n_\beta = \lim _{K \rightarrow \infty } \sum
      _{|m| \leq K} \frac{1}{2\pi } \int _{A_m} \omega ^n_\beta = 0 .
   \]
\end{Coro}

Corollary~\ref{Corollary 6.1} can be combined with
Proposition~\ref{Proposition 5.4} yielding the following

\begin{Prop}
\label{Proposition 6.2} Let $n \in {\mathbb Z}$ and $\beta = (\beta _k)_{k
\not= n} \in \ell ^1_{\hat n}$ be arbitrary. If $\int _{A_m} \omega ^n_\beta
= 0$ for any $m \in {\mathbb Z} \backslash \{ n\}$, then $\beta = 0$.
\end{Prop}


\section{Existence of normalized differentials}
\label{7. Existence of normalized differentials}

The aim of this section is to study the operators $T^n$ and $T^\ast $, introduced
in Section~\ref{3. Outline of proof of Theorem 1.2}, and to prove
Proposition~\ref{Proposition 4.2}, Proposition \ref{Proposition 4.3}, and 
Proposition~\ref{Proposition 4.4}. We begin with the study of $T^\ast $.

\begin{Lemma}\hspace{-2mm}{\bf .}
\label{Lemma 7.1} Locally uniformly on ${\mathcal W}$, the coefficients $T^\ast
_{mj}$ of $T^\ast $ satisfy the following estimates
   \[ T^\ast _{mj} = \begin{cases} O \Big( \frac{|\dot \lambda _m - \tau _m| +
      |\gamma _m|}{j - m} \Big), & j \in\Z,\, |m|> N_0,\, m \not= j \\ 
     1 + \ell ^2(m) , &|j| > N_0,\, m = j \\ 
     O \Big( \frac{1}{1 + |j|} \Big), &j \in {\mathbb Z},\, |m| \leq N_0 
               \end{cases}
   \]
\end{Lemma}

\noindent{\em Proof.} {\em Case $|j| > N_0$ and $\lambda ^+ _m = \lambda ^-_m$:}
Recall that by Proposition~\ref{Proposition 2.2} and the definition of the
canonical root
   \[ \frac{\dot \Delta (\lambda )}{\sqrt[c]{{\mathcal R}(\lambda )}} = \frac{1}
      {i} \prod _{k \in {\mathbb Z}} \frac{\dot \lambda _k - \lambda }{\sqrt[s]
      {(\lambda ^+_k - \lambda )(\lambda ^-_k - \lambda )}}
   \]
and that for $|j| > N_0$, $T^\ast _{mj}$ is given by
   \[ T^\ast _{mj} = \frac{1}{2\pi } \int _{\Gamma _m} \frac{1}{\lambda - \dot
      \lambda _j} \frac{\dot \Delta (\lambda )}{\sqrt[c]{{\mathcal R}(\lambda )}}\,
      d \lambda .
   \]
As $\lambda ^+_m = \lambda ^-_m = \tau _m=\dot\lambda_m$ one has by the definition of the
standard root, $\sqrt[s]{(\lambda ^+_m - \lambda )(\lambda ^-_m - \lambda )}= \tau _m - \lambda $, 
and hence, if in addition $m \not= j$,
$\frac{1}{\lambda - \dot \lambda _j} \frac{\dot \Delta (\lambda )}{\sqrt[c]{{\mathcal
R}(\lambda )}}$ is holomorphic near $\tau _m$ and thus
$T^\ast _{mj} = 0$. If $m = j$ (and hence $|m| > N_0$) one gets
\begin{eqnarray} 
T^\ast _{mm}&=& \frac{1}{2\pi} \int _{\Gamma _m} \frac{1}{\lambda -
                   \dot \lambda _m} \frac{\dot\Delta(\lambda )}{\sqrt[c]{{\mathcal R}(\lambda )}}\,d\lambda\\
            &=& \frac{1}{2\pi i} \int _{\Gamma _m} \frac{1}{\lambda-\tau_m}\prod _{k \not= m}
                \frac {\dot \lambda _k - \lambda }{\sqrt[s]{(\lambda ^+_k - \lambda )
                  (\lambda ^-_k - \lambda )}}\,d\lambda .
\end{eqnarray}
Therefore, by the residue theorem and the product estimate in
\cite[Lemma C.3]{GKP}
   \[ T^\ast _{mm} = \prod _{k\not= m} \frac{\dot \lambda _k - \tau _m}{\sqrt[s]
      {(\lambda ^+_k - \tau _m)(\lambda ^-_k - \tau _m)}} = 1 + \ell ^2 (m)
   \]
locally uniformly in $\W$.

\medskip

\noindent {\em Case $|j| > N_0$ and $\lambda ^+_m \not= \lambda ^-_m$:} If $m = j$, 
one uses again \cite[Lemma C.3]{GKP} to see that
   \begin{align*} T^\ast _{mm} &= \frac{1}{2\pi i} \int _{\Gamma _m} \frac{-1}
                     {\sqrt[s]{(\lambda ^+_m - \lambda )(\lambda ^-_m -
                     \lambda )}} \prod _{k\not= m} \frac{\dot \lambda _k -
                     \lambda }{\sqrt[s]{(\lambda ^+_k - \lambda )(\lambda ^-_k
                     - \lambda )}}\,d\lambda \\
                  &= -\frac{1}{2\pi i} \int _{\Gamma _m} \frac{d\lambda}{\sqrt[s]
                     {(\lambda ^+_m - \lambda )(\lambda ^-_m - \lambda )}} +
                     \ell ^2(m) .
   \end{align*}
A direct calculation shows that
$\frac{1}{2\pi i} \int _{\Gamma _m} \frac{d\lambda}{\sqrt[s]{(\lambda ^+_m - \lambda )
(\lambda ^-_m - \lambda )}} = -1$
yielding the claimed estimate $T^\ast_{mm} = 1 + \ell ^2(m)$ in this case. 
If $|m| > N_0$ but $m \not= j$, then
\begin{equation}\label{eq:T^*_mj}
T^\ast _{mj} = \frac{1}{2\pi i}\int_{\Gamma_m}\frac{B_{mj}(\lambda)\,d\lambda}
{\sqrt[s]{(\lambda ^+_m - \lambda )(\lambda ^-_m - \lambda )}}
\end{equation}
where
\begin{equation}\label{eq:B_mj}
B_{mj}:=-\frac{\dot\lambda_m - \lambda }{\sqrt[s]{(\lambda^+_j - \lambda )(\lambda^-_j-\lambda )}}
      \prod_{k\not= m,j} \frac{\dot \lambda _k
      - \lambda }{\sqrt[s]{(\lambda ^+_k - \lambda )(\lambda ^-_k - \lambda )}}.
\end{equation}
By deforming the contour $\Gamma_m$ to the straight interval $G_m$ (taken twice) one then sees
(cf. \cite[Lemma 14.3]{GKP}) that
   \[ |T^\ast _{mj}| \leq \max _{\lambda \in G_m}|B_{mj}(\lambda)|.
   \]
For $\lambda \in G_m$, 
$|\dot \lambda _m - \lambda | \leq |\dot \lambda _m -\tau _m| + |\gamma _m|$ and 
\[
\Big( \sqrt[s]{(\lambda ^+_j - \lambda )(\lambda ^-_j- \lambda ) }\Big) ^{-1} =O\Big(\frac{1}{j - m}\Big)
\]
whereas again with \cite[Lemma C.3]{GKP}, 
$\prod\limits_{k \not= m,j}\frac{\dot\lambda_k-\lambda}{\sqrt[s]{(\lambda^+_k-\lambda)(\lambda^-_k-\lambda )}}=O(1)$
uniformly on $\W$. Altogether we thus conclude that in the case considered
$T^\ast _{mj} = O \Big( \frac{|\dot \lambda _m - \tau _m| + |\gamma _m|}{j-m}\Big)$. 
Finally if $|m| \leq N_0$ and hence $m \not= j$ (as we assume $|j| >N_0$), one has that
\eqref{eq:T^*_mj} and \eqref{eq:B_mj} hold. 
Hence,
\[ 
|T^\ast _{mj}|\le 
\max_{\lambda\in\Gamma_m}\left|\frac{B_{mj}(\lambda)}
{\sqrt[s]{(\lambda ^+_m - \lambda )(\lambda ^-_m - \lambda )}}\right|\,\mathop{\rm length}(\Gamma_m)/2\pi
\]
where $\mathop{\rm length}(\Gamma_m)$ is the Euclidean length of $\Gamma_m$.
Using that $\lambda ^\pm _m$ are inside $\Gamma_m'$ and hence uniformly in
${\mathcal W}$ separated from $\Gamma _m$ and as by definition, different
contours are apart by a uniform constant one concludes from the estimate
$\Big( \sqrt[s]{(\lambda ^+_j-\lambda )(\lambda^-_j-\lambda)}\Big)^{-1}=O\Big(\frac{1}{j}\Big)$ and 
\cite[Lemma C.3]{KP} that $T^\ast_{mj}= O\big(\frac{1}{j}\big)$ uniformly on $\W$.

\medskip

\noindent{\em Case $|j| \leq N_0$:} In this case the coefficient $T^\ast _{mj}$
is given by the formula
\begin{align*} T^\ast _{mj} &= \frac{1}{2\pi } \int _{\Gamma _m} \lambda^{N_0
                     + j} \frac{\dot \Delta (\lambda )}{\prod\limits_{|k| \leq N_0}
                     (\lambda - \dot \lambda _k)} \frac{d\lambda}{\sqrt[c]{
                     {\mathcal R}(\lambda )}} \\
                  &= \frac{1}{2\pi i} \int_{\Gamma _m}\!\!\!\frac{- \lambda ^{N_0
                     + j}}{\prod\limits_{|k|\leq N_0}\!\!\!\sqrt[s]{(\lambda ^+_k -
                     \lambda )(\lambda ^-_k - \lambda )}} \prod _{|k| > N_0}
                     \frac{\dot \lambda _k - \lambda }{\sqrt[s]{(\lambda ^+_k
                     - \lambda )(\lambda ^-_k - \lambda )}}\,d\lambda\,.
\end{align*}
If $|m|>N_0$ we apply \cite[Lemma 14.3]{KP} to conclude that,
\[ 
|T^*_{mj}|\le\!\!\max_{\lambda\in G_m}\!\!\Big|  
              \frac{-\lambda ^{N_0+ j}(\dot\lambda_m-\lambda)}{\prod\limits_{|k|\leq N_0}\!\!\!
                \sqrt[s]{(\lambda ^+_k -\lambda )(\lambda ^-_k - \lambda )}} \!\!
                  \prod_{|k| > N_0,k\ne m}\!\!
                     \frac{\dot \lambda _k - \lambda }{\sqrt[s]{(\lambda ^+_k
                     - \lambda )(\lambda ^-_k - \lambda )}}\Big|.
\]
As for $|m| > N_0$, one has
   \[ \max _{\lambda \in G _m }\Big\arrowvert \prod _{|k| \leq N_0}
       \frac{\lambda ^{N_0 + j}}
      {\sqrt[s]{(\lambda ^+_k - \lambda )(\lambda ^-_k - \lambda )}}
      \Big\arrowvert = O \Big( \frac{1}{m} \Big)\,,
   \]
it follows by the product estimate \cite[Lemma C.3]{GKP} that 
\[
T^\ast_{mj}=O\Big(\frac{|\dot\lambda_m-\tau_m|+|\gamma_m|}{m}\Big)
\]
uniformly on $\W$. If $|m| \leq N_0$ we use again that
$\lambda ^\pm _m$ are inside $\Gamma '_m$, hence uniformly in ${\mathcal W}$
separated form $\Gamma _m$, and that different contours are apart by a
uniform constant, to see that $T^\ast _{mj} = O(1)$. 
The claimed estimates for $T^*_{mj}$ are proved.
\finishproof

\noindent From Lemma~\ref{Lemma 7.1} it immediately follows that $T^\ast $ defines a
bounded linear operator, $T^\ast : \ell ^1 \rightarrow \ell ^1$.

\smallskip

\noindent{\em Proof of Proposition~\ref{Proposition 4.4}. } Take $\varphi\in\W$.
By Proposition~\ref{Proposition 5.1}, $T^*$ is injective. We claim that $T^*-Id$
is a compact operator on $\ell^1$. Therefore, $T^\ast $ is Fredholm and thus
$T^\ast$ is a linear isomorphism. To see that $T^\ast-Id$ is compact
introduce for any $N > N_0$ the operator $K_N : \ell ^1 \rightarrow \ell ^1$, 
\[
K_N :=\Pi_N\circ\big(T^\ast - Id\big),
\]
where $\Pi _N : \ell ^1\rightarrow \ell ^1$ is the projection,
\[
(\beta _k)_{k \in {\mathbb Z}}\mapsto (\cdots , 0, \beta _{-N},\cdots,\beta _N, 0,\cdots ).
\] 
Note that $K_N$ is an operator of finite rank and therefore compact. 
By Lemma~\ref{Lemma 7.1} we have for any $N \geq N_0$,
   \[ 
(T^\ast - Id - K_N)_{mj} = 
\begin{cases} 
        O\Big( \frac{|\dot \lambda _m - \tau _m| + |\gamma_m|}{j-m} \Big), 
             &j\in{\mathbb Z} \backslash \{ m\},\,|m| > N \\ 
        \ell ^2(m), &|m| > N,\,j = m \\
           0, &j \in{\mathbb Z},\,|m| \leq N\,.
\end{cases}
   \]
Hence there exists a constant $C > 0$ independent of $N \geq N_0$ so that
for any $N \geq N_0$ and any $\beta \in \ell ^1$,
\begin{eqnarray*} \| (T^\ast - Id - K_N) \beta \| _{\ell ^1} &\le& C \sum _{|m| > N}
      \sum _{j \not= m} \frac{|\dot \lambda _m - \tau _ m| + |\gamma _m|}
      {|j-m|} |\beta _j| \\
      &+&\sum_{|m|>N} \ell^2(m) |\beta_m|.
\end{eqnarray*}
Clearly,
\[
\sum_{|m|>N}\ell^2(m) |\beta_m|=\Big(\sup_{|m|>N}|\ell^2(m)|\Big) \|\beta\|_{\ell^1}\to 0
\]
as $N\to\infty$. By changing the order of summation we get from the Cauchy-Schwartz inequality
\[ \sum _{|m| > N} \sum _{j \not= m} \frac{|\dot \lambda _m - \tau _m|
      + |\gamma _m|}{|j - m|} |\beta _j| \leq C \| \beta \| _{\ell ^1}
      \big(\!\!\!\sum _{|m| > N} |\dot \lambda _m - \tau _m|^2 + |\gamma _m|^2
      \big)^{1/2} .
\]
Altogether it then follows that
   \[ \| T^\ast - Id - K_N \|_{{\mathcal L}(\ell ^1)} \rightarrow 0 \,\,\mbox { as } N\rightarrow\infty,
   \]
showing that $T^\ast - Id$ is compact.
\finishproof

Next we want to prove Proposition~\ref{Proposition 4.2}. First we
establish the following two lemmas

\begin{Lemma}\label{Lemma 7.2} 
For any $n \in {\mathbb Z}$ and for any $\varphi\in{\mathcal W}$ 
the coefficients $b^n_m$, $m \not= n$, of $b^n$ satisfy
   \[ b^n_m-b^*_m = \begin{cases} O \Big( \frac{|\dot \lambda _m - \tau _m| +
      |\gamma _m|}{m - n}\Big), &|m| > N_0,\, \ m \not= n \\ 
       O\Big(\frac{1}{n}\Big), &|m| \leq N_0,\,m \not= n  \end{cases}
   \]
uniformly in $n\in\Z$ and locally uniformly in $\W$.  
In particular, it follows that $b^n$ is in $\ell ^1_{\hat n}$.
\end{Lemma}

\noindent{\em Proof.} As the cases $|n| > N_0$ and $|n| \leq N_0$ can be
treated in the same way, we only consider the case $|n| > N_0$.

\medskip

\noindent {\em Case $|m| > N_0$:} Taking into account that 
$b^*_m=\frac{1}{2\pi}\int_{\Gamma_m} \frac{\dot \Delta (\lambda )}{\sqrt[c]{{\mathcal R}(\lambda )}}\,d\lambda$
and 
$b^n_m = \frac{1}{2\pi }\int _{\Gamma _m} \frac{m\pi - n\pi }{\lambda - \dot \lambda _n} \frac{\dot
\Delta (\lambda )}{\sqrt[c]{{\mathcal R}(\lambda )}}\,d\lambda $
we see that
   \[ b^n_m-b^*_m = \frac{1}{2\pi i} \int _{\Gamma _m} \Big( \frac{m \pi - n\pi }
      {\lambda - \dot \lambda _n} - 1 \Big) \frac{\dot \lambda _m - \lambda }
      {\sqrt[s]{(\lambda ^+_m - \lambda )(\lambda ^-_m - \lambda )}}\,\Pi_m(\lambda)\,d\lambda
   \]
where $\Pi_m(\lambda):=\prod\limits_{k \not= m} \frac{\dot \lambda _k - \lambda }{\sqrt[s]{(\lambda ^+_k -
      \lambda )(\lambda ^-_k - \lambda )}}$.
If $\lambda ^+_m = \lambda ^-_m$, one has $\dot\lambda _m = \lambda^+_m=\lambda^-_m$ and
$\sqrt[s]{(\lambda ^+_m - \lambda )(\lambda ^-_m - \lambda )} = \dot\lambda_m -\lambda $. 
As $m \not= n$ one then concludes from the analyticity of the integrand in 
$D_m(\pi/4)$ that $b^n_m = 0$. If $\lambda ^+_m \not= \lambda ^-_m$ we first note that
   \[ \frac{m\pi - n\pi }{\lambda - \dot \lambda _n} - 1 = \frac{(m\pi -
      \lambda ) + (\dot \lambda _n - n \pi )}{\lambda - \dot \lambda _n}
      = O\Big(\frac{1}{m - n}\Big)
   \]
where we used that $\dot \lambda _n - n\pi = O(1)$ by
Proposition~\ref{Proposition 2.1}. Furthermore, by \cite[Lemma C.3]{GKP}, $\Pi_m(\lambda) = O(1)$. 
By deforming the contour $\Gamma _m$ to the straight interval $G_m$ (taken twice) and by using that
for $\lambda \in G_m$, $|\dot \lambda _m - \lambda | \leq |\dot \lambda _m - \tau _m| + |\gamma _m|$,
one sees from \cite[Lemma 14.3]{KP} that 
$b^n_m-b^*_m = O \big( \frac{|\dot \lambda _m - \tau _m| + |\gamma _m|}{m - n} \big)$.

\medskip

\noindent {\em Case $|m| \leq N_0$:} We can argue similarly as
above to conclude that $b^n_m-b^*_m = O \big( \frac{1}{n}\big)$ (cf. Lemma \ref{Lemma 7.1}). 

Going through the arguments of the proof one verifies that the estimates obtained
are uniform in $n\in\Z$ and locally uniform on ${\mathcal W}$.
\finishproof

The coefficients $T^n_{mj}$ can be estimated using Lemma~\ref{Lemma 7.1} by
writing $T^n_{mj} = T^\ast _{mj} + R^n_{mj}$ $(m,j\ne n)$ where $R^n_{mj}$ is defined
for $|n| > N_0$ as follows
   \[ R^n_{mj} := \begin{cases} \frac{1}{2\pi } \int _{\Gamma _m} \Big(
      \frac{m\pi - n \pi }{\lambda - \dot \lambda _n} - 1 \Big) 
      \frac{\lambda^{N_0 + j}}{\prod\limits_{|k| \leq N_0}(\lambda - \dot\lambda _k)}
      \frac{\dot
      \Delta (\lambda )}{\sqrt[c]{{\mathcal R}(\lambda )}}\,d\lambda , 
      &|j| \leq N_0 \\ 
       \frac{1}{2\pi } \int _{\Gamma _m} \Big(
      \frac{m\pi - n \pi }{\lambda - \dot \lambda _n} - 1 \Big) \frac{\dot
      \Delta (\lambda )}{\lambda - \dot \lambda _j}
      \frac{d\lambda }{\sqrt[c]{{\mathcal R}(\lambda )}} , &
      |j| > N_0 \end{cases}
   \]
whereas in the case $|n| \leq N_0$,\footnote{See \eqref{eq:epsilon} for the definition of $\varepsilon ^n_j$.}
   \[ R^n_{mj} := \begin{cases} \frac{1}{2\pi } \int _{\Gamma _m} \Big(
      \frac{m\pi - N_0 \pi }{\lambda - \dot \lambda _{N_0}} - 1 \Big) 
      \frac{\lambda ^{N_0 + j - \varepsilon ^n_j}}
      {\prod\limits_{|k| \leq N_0, k\not= N_0}
      (\lambda - \dot \lambda _k)}
      \frac{\dot\Delta (\lambda )}{\sqrt[c]{{\mathcal R}
      (\lambda )}}\, d\lambda , &|j| \leq N_0 \\
      \frac{1}{2\pi } \int _{\Gamma _m} \Big(
      \frac{m\pi - N_0 \pi }{\lambda - \dot \lambda _{N_0}} - 1 \Big) \frac{\dot
      \Delta (\lambda )}{\lambda - \dot \lambda _j}
      \frac{d\lambda }{\sqrt[c]{{\mathcal R}(\lambda )}} , &|j| > N_0 .\end{cases}
   \]
Note that by Proposition~\ref{Proposition 2.2}, for $\lambda \in \Gamma _m$ and $|n|>N_0$,
   \[ \frac{m\pi - n\pi }{\lambda - \dot \lambda _n} - 1 = \frac{(m \pi -
      \lambda ) + (\dot \lambda _n - n \pi )}{\lambda - \dot \lambda _n}
      = O \Big( \frac{1}{m-n} \Big) .
   \]
It is convenient to rewrite $R^n_{mj}$ in the case $|n| > N_0$ as follows
\begin{equation} \label{eq:****}
      \begin{cases} \frac{1}{2\pi i} \int _{\Gamma _m}
      \frac{(m\pi - \lambda ) + (\dot \lambda _n -n\pi  )}{\sqrt[s]{
      (\lambda ^+_n - \lambda )(\lambda ^-_n - \lambda)}} \frac{
      \lambda ^{N_0 + j}}{\underset{|k| \leq N_0}{\prod } \sqrt[s]{(\lambda ^+_k -
      \lambda )(\lambda ^-_k - \lambda )}}\,
      \Pi_n(\lambda)\,d\lambda , &|j|\leq N_0 \\ 
      \frac{1}{2\pi i} \int _{\Gamma _m}\frac{(m\pi - \lambda ) + (\dot \lambda _n -n\pi  )}{\sqrt[s]{
      (\lambda ^+_n - \lambda )(\lambda ^-_n - \lambda)}} \frac{1}
      {\sqrt[s]{(\lambda ^+_j - \lambda )(\lambda ^-_j - \lambda )}}\,
      \Pi_{nj}(\la)\,d\lambda , &|j|> N_0  
      \end{cases}
\end{equation}
where 
\[
\Pi_n(\lambda):=\prod\limits_{|k|> N_0, k\ne n}
\frac{\dot\lambda _k - \lambda }{\sqrt[s]{(\lambda ^+_k -\lambda )(\lambda ^-_k - \lambda )}}
\]
and
\begin{equation}\label{eq:Pi_nj}
\Pi_{nj}(\lambda):=\prod\limits_{k\ne n,j}\frac{\dot \lambda _k - \lambda }{\sqrt[s]{(\lambda ^+_k -
      \lambda )(\lambda ^-_k - \lambda )}}.
\end{equation}
In the case $|n| \leq N_0$, similar identities hold.

\begin{Lemma}\hspace{-2mm}{\bf .}\label{Lemma 7.3} 
For any $n\in\Z$ and for any $\varphi\in\W$, the coefficients $R^n_{mj}$ satisfy the following estimates
   \[ R^n_{mj} = \begin{cases} O \Big( \frac{|\dot \lambda _m -\tau _m| + |\gamma
      _m|}{(m - n)(j - m)} \Big), &j \in {\mathbb Z}\setminus\{ n \},\,|m| > N_0,\,m \ne j, n \\ 
      \frac{\ell ^2(m) + \ell ^2(n)}{m - n}, &|j| > N_0,\, \ m = j,\,m \ne n \\ 
       O\Big(\frac{1}{(m - n)(1 + |j|)}\Big), &j \in {\mathbb Z}\setminus\{ n \},\,
      |m| \leq N_0,\,m \not= n\end{cases}
   \]
uniformly in $n$ and locally uniformly on $\W$.
\end{Lemma}

\noindent{\em Proof.} As the cases $|n| > N_0$ and $|n| \leq N_0$ are proved in the
same way we will only consider the case $|n| > N_0$.

Throughout the proof we assume that $m \ne n$ and $j \not= n$. The proof is very
similar to the one of Lemma~\ref{Lemma 7.1}.

\noindent {\em Case $|j| > N_0$ and $\lambda ^+_m = \lambda ^-_m (=\tau _m)$:}
We argue as in the proof of Lemma~\ref{Lemma 7.1} to see that for $m \not= j$, $R^n_{mj} = 0$. 
If $m = j$ (and hence $|m| > N_0$) one gets from the residue theorem and \cite[Lemma C.3]{GKP}
that
   \begin{align*} R^n_{mm} &= \frac{(m\pi - \tau _m) + (\dot \lambda _n - n \pi)}
                     {\sqrt[s]{(\lambda ^+_n - \tau _m)(\lambda ^-_n - \tau _m)}}
                     \prod _{k\not= m,n} \frac{\dot \lambda _k - \tau _m}{\sqrt
                     [s]{(\lambda ^+_k - \tau _m)(\lambda ^-_k - \tau _m)}} \\
                  &= O \Big( \frac{(m\pi - \tau_m) + (\dot \lambda _n - n\pi )}
                     {m - n} \Big) .
   \end{align*}
By Proposition~\ref{Proposition 2.2}, $m \pi - \tau _m = \ell ^2(m)$ and
$\dot \lambda _n - n\pi = \ell ^2 (n)$, hence $R^n_{mm} = \frac{\ell ^2(m) +
\ell ^2(n)}{m - n}$.

\medskip

\noindent {\em Case $|j| > N_0$ and $\lambda ^+_m \not= \lambda ^-_m$:} If $j = m$ 
(and hence $|m| > N_0$) one deforms the contour $\Gamma _m$ to the straight interval $G_m$ (taken twice)
and obtains from \cite[Lemma 14.3]{GKP} that
   \begin{align*} R^n_{mm} &= \frac{1}{2\pi i} \int_{\Gamma _m}\!\!\!\frac{(m\pi -
                     \lambda ) + (\dot \lambda _n - n\pi )}{\sqrt[s]{(\lambda ^+
                     _n - \lambda )(\lambda ^-_n - \lambda )}} \frac{1}{\sqrt[s]
                     {(\lambda ^+_m - \lambda )(\lambda ^-_m - \lambda )}}\,
                     \Pi_{nm}(\lambda)\,d\lambda \\
                      &= \frac{\ell ^2(m) + \ell ^2(n)}{m - n}
   \end{align*}
where $\Pi_{mn}(\lambda)$ is defined in \eqref{eq:Pi_nj}.
If $|m| > N_0$, but $m \not= j$, one argues similarly: Deforming the contour
$\Gamma_m$ to $G_m$ (taken twice) one sees that 
\begin{equation*} 
|R^n_{mj}|\le\max_{\lambda\in G_m}\left|\frac{(m\pi-\lambda )
                     + (\dot \lambda _n - n \pi )}{\sqrt[s]{(\lambda ^+_n - \lambda )
                     (\lambda ^-_n - \lambda )}}\,
                     \frac{({\dot\lambda}_m-\lambda)}{\sqrt[s]{(\lambda_j^+-\lambda)(\lambda_j^--\lambda)}}\,
                     \Pi_{mjn}(\lambda)\right|
\end{equation*}
where $\Pi_{mjn}(\lambda):=\prod\limits_{k\ne m,j,n}\frac{\dot \lambda _k - \lambda }{\sqrt[s]{(\lambda ^+_k -
      \lambda )(\lambda ^-_k - \lambda )}}$.
As for $|m| > N_0$ and $\lambda\in G_m$, 
$|\dot \lambda _m - \lambda | \leq|\dot \lambda _m - \tau _m| + |\gamma _m |$ and 
$\Big(\sqrt[s]{(\lambda ^+_j - \lambda )(\lambda ^-_j - \lambda )}\Big)^{-1}= O\Big(\frac{1}{j-m}\Big)$,
one concludes that 
$R^n_{mj} = O \Big( \frac{|\dot \lambda_m - \tau _m| + |\gamma _m|}{(m - n)(j-m)} \Big)$. 
Finally, if $|m| \leq N_0$ and hence $m \not= j$ (as we assume $|j| > N_0$), one argues as in the proof
of Lemma~\ref{Lemma 7.1} to conclude that
   \[ R^n_{mj} = O \Big(\frac{1}{(m-n)(1+|j|)}\Big)\,.
   \]

\medskip

\noindent {\em Case $|j| \leq N_0$:} In this case $R^n_{mj}$ is given by the first equation 
in \eqref{eq:****}. If  $|m| > N_0$ note that
   \[ \max _{\lambda \in G_m} \left| \frac{\lambda ^{N_0 + j}}
      {\prod _{|k| \leq N_0} \sqrt[s]{(\lambda ^+_k - \lambda )(\lambda ^-_k -
      \lambda )}}\right| = O \Big( \frac{1}{m} \Big)\,,
   \]
   \[ \prod _{|k| > N_0, k\ne n,m} \frac{\dot \lambda _k - \lambda}
      {\sqrt[s]{(\lambda ^+_k - \lambda )(\lambda ^-_k - \lambda )}} = O (1)\,,
   \]
and $|\dot \lambda _m - \lambda| \leq |\dot \lambda _m - \tau _m | + |\gamma _m|$
for any $\lambda \in G_m$. This together with \cite[Lemma 14.3]{GKP} implies that,
   \[ R^n_{mj} = O \Big( \frac{|\dot \lambda _m - \tau _m| + |\gamma _m|}{(m-n)
      m} \Big)\left( = O \Big( \frac{|\dot \lambda _m - \tau _m| + |\gamma _m|}
      {(m - n)(j - m)} \Big) \right) .
   \]
In the remaining case $|m|\le N_0$ we argue again as in the proof of
Lemma~\ref{Lemma 7.1} to see that $R^n_{mj} = O(1)$.
The claimed estimates for $R^n_{mj}$ are thus proved. Going through the
arguments of the proofs one sees that the derived estimates hold uniformly in $n$ and locally uniformly
on ${\mathcal W}$.
\finishproof

\medskip

\noindent From Lemma~\ref{Lemma 7.1} and Lemma~\ref{Lemma 7.3} it immediately
follows that for each $n \in {\mathbb Z}$, $\ T^n$ defines a bounded linear
operator, $T^n : \ell ^1_{\hat n} \rightarrow \ell ^1_{\hat n}$.

\medskip

\noindent{\em Proof of Proposition~\ref{Proposition 4.2}.} By
Proposition~\ref{Proposition 6.2}, $T^n$ is injective. Arguing as in the proof of
Proposition~\ref{Proposition 4.4} one sees that $T^n - Id$ is a compact
operator on $\ell^1_{\hat n}$. Therefore $T^n$ is Fredholm and thus $T^n$
a linear isomorphism. Finally, by Lemma~\ref{Lemma 7.2}, $b^n \in \ell ^1
_{\hat n}$.
\finishproof

\noindent Now let us turn to the proof of Proposition~\ref{Proposition 4.3}.
Recall that for any $n \in {\mathbb Z}$ and for any given $\varphi$ in
${\mathcal W}$ we denote by $\beta ^n\equiv\beta^n(\varphi)$ the unique solution of 
$T^n \beta = b^n$. In this way we obtain maps 
\[
\beta ^n : {\mathcal W} \to\ell ^1_{\hat n}
\] 
and 
\[
\zeta _n : {\mathbb C} \times {\mathcal W}\rightarrow {\mathbb C},\,\,\,\,
(\lambda , \varphi ) \mapsto \zeta ^n_{\beta^n}(\lambda )\,.
\] 
Similarly,  for any given $\varphi\in\W$ we denote by $\beta^*\equiv\beta^*(\varphi)$
the unique solution of the linear system $T^*\beta=b^*$ (see Proposition \ref{Proposition 4.4}).
Recall that by Lemma \ref{lem:zero_form}, $b^*$ is in $l^1$ and independent of $\varphi\in\W$.
In this way we obtain a map
\[
\beta^* : \W\to\ell^1\,.
\]
\begin{Lemma}\hspace{-2mm}{\bf .}\label{Lemma 7.4} 
For any $n \in {\mathbb Z}$, (i) $b^n : {\mathcal W}
\rightarrow \ell ^1_{\hat n}$ and (ii) $T^n : {\mathcal W} \rightarrow {\mathcal L}(\ell ^1_{\hat n})$ 
are analytic maps. In addition, (iii) $T^* : \W\rightarrow {\mathcal L}(\ell^1)$ is an analytic map.
\end{Lemma}
Combining Lemma \ref{Lemma 7.4} with Proposition \ref{Proposition 4.2} and 
Proposition \ref{Proposition 4.4} one gets
\begin{Coro}\label{coro:beta_analytic}
For any $n\in\Z$ the maps $\beta^n : \W\to\ell ^1_{\hat n}$ and
$\beta^* : \W\to\ell^1$ are analytic.
\end{Coro}
\noindent{\em Proof of Lemma \ref{Lemma 7.4}.} 
(i) Let us first consider the case $|n| > N_0$. According to
the definition of $b^n = (b^n_m)_{m\not= n}$ in Section~\ref{3. Outline of proof of Theorem 1.2},
   \[ b^n_m = \frac{1}{2\pi } \int _{\Gamma _m} \frac{m\pi - n\pi}{\lambda -
      \dot \lambda _n} \frac{\dot \Delta (\lambda )}{\sqrt[c]{\Delta ^2
      (\lambda ) - 4}}\,d\lambda .
   \]
By \cite[Theorem A.3]{GKP} it suffices to show that $b^n$ is locally bounded
and weakly analytic. By Lemma~\ref{Lemma 7.2} and the Cauchy-Schwartz estimate,
$b^n$ is locally bounded. As the dual of $\ell ^1_{\hat n}$ is 
$\ell ^\infty _{\hat n} \equiv \ell ^\infty({\mathbb Z} \backslash \{ n \},{\mathbb C})$, 
in view of Montel's theorem, the weak analyticity of $b^n$ then follows
once we prove that each component $b^n_m$, $m \not= n$, of $b^n$ is analytic on ${\mathcal W}$. 
To this end let us analyze the integrand in the definition of $b^n_m$. Recall that
$\Delta , \dot \Delta : {\mathbb C} \times L^2_c \rightarrow {\mathbb C}$
are analytic maps. As $|n| > N_0$ by assumption, $\dot \lambda _n$ is a
simple zero of $\dot \Delta (\lambda )$ (cf. Proposition~\ref{Proposition 2.1}) 
and hence we obtain from the implicit function theorem that
$\dot \lambda _n : {\mathcal W} \rightarrow {\mathbb C}$ is analytic. 
By construction, $\sqrt[c]{\Delta ^2(\lambda ) - 4}$ is analytic on 
$G\backslash \sqcup _{k \in {\mathbb Z}} G_k$. 
In view of the definition of $\Gamma_m$ and $G_m$, $m\in{\mathbb Z}$ (Section~\ref{2. Preliminaries})
there exists $\varepsilon > 0$, independent of $m$, so that $\sqrt[c]{\Delta^2 - 4}$ is analytic on 
$U_\varepsilon (\Gamma _m) \times {\mathcal W}$,
where $U_\varepsilon (\Gamma _m)$ is the $\varepsilon $-tubular neighborhood
of $\Gamma _m$, 
$U_\varepsilon (\Gamma _m) := \{\lambda\in{\mathbb C}\,|\,\mathop{\rm dist}(\lambda,\Gamma_m)<\varepsilon\}$. 
It then follows that for any $m \not= n$, $b^n_m : {\mathcal W} \rightarrow {\mathbb C}$ is
analytic.

\medskip

\noindent In the case $|n| \leq N_0$, $b^n_m$ is defined for any $m \not= n$ by
   \[ b^n_m = \frac{1}{2\pi } \int _{\Gamma _m} \frac{m \pi - N_0 \pi }{\lambda
      - \dot \lambda _{N_0}} \frac{\dot \Delta (\lambda )}{\sqrt[c]{\Delta ^2
      (\lambda ) - 4}}\,d\lambda .
   \]
By the choice of $N_0$ in Proposition~\ref{Proposition 2.1}, for any
$\varphi \in {\mathcal W}, \dot \lambda _{N_0}$ is the only root of $\dot \Delta
(\cdot , \varphi )$ in $D_{N_0}(\pi / 4)$. Hence arguing as above, $\dot
\lambda _{N_0} : {\mathcal W} \rightarrow {\mathbb C}$ is analytic. Using
the same arguments as in the case $|n| > N_0$ one sees that $b^n : {\mathcal W}
\rightarrow {\mathbb C}$ is analytic also in this case.

\medskip

\noindent (ii) As the proofs for $|n| > N_0$ and $|n| \leq N_0$ are similar,
we consider the case $|n| > N_0$ only. In this case, the coefficients
$T^n_{mj}$ $(m, j \in {\mathbb Z} \backslash \{ n \})$ of $T^n$ are given by
\[ 
T^n_{mj} = 
\begin{cases}
\frac{1}{2\pi }\int_{\Gamma _m} 
\frac{m\pi-n\pi}{\lambda-\dot\lambda_n}\frac{1}{\lambda -\dot\lambda_j} 
\frac{\dot \Delta (\lambda )}{\sqrt[c]{\Delta^2(\lambda)-4}}\,d\lambda,&
|j| > N_0 \\
\frac{1}{2\pi}\int_{\Gamma _m}\frac{m\pi-n\pi}{\lambda-\dot\lambda_n} 
\frac{\dot\Delta(\lambda)}{\prod_{|k|\le N_0}(\lambda-\dot\lambda_k)}
\frac{\lambda^{N_0+j}}{\sqrt[c]{\Delta^2(\lambda )-4}}\,d\lambda,&
|j|\le N_0\,.
\end{cases}
\]
Again by \cite[Theorem A.3]{GKP} it suffices to show that $T^n$ is locally
bounded and weakly analytic. By Lemma~\ref{Lemma 7.1} and Lemma~\ref{Lemma
7.3}, $T^n$ is locally bounded (see the arguments in the proof of 
Proposition \ref{Proposition 4.4}). It remains to show that $T^n$ is weakly
analytic. First note that by arguing as in (i) one sees that for any $m,
j \in {\mathbb Z} \backslash \{ n\}, T^n_{mj}: {\mathcal W} \rightarrow
{\mathbb C}$ is analytic. Since $T^n : {\mathcal W} \rightarrow
{\mathcal L}(\ell ^1_{\hat n})$ is locally bounded so is for any $m \not=
n$ the map
   \[ T^n_m : {\mathcal W} \rightarrow \ell ^\infty _{\hat n}, \varphi
      \mapsto (T^n_{mj}(\varphi ))_{j \not= m} .
   \]
Using that the components $T^n_{mj}$ of $T^n_m$ are analytic it then
follows from \cite[Theorem A.3]{KP} that $T^n_m : {\mathcal W} \rightarrow
\ell ^\infty _{\hat n}$ is analytic. As a consequence, for each $N > |n|$,
the map $\Pi _N \circ T^n : {\mathcal W} \rightarrow {\mathcal L}(\ell ^1_{\hat n})$ 
is analytic, where $\Pi _N : \ell ^1_{\hat n}\rightarrow \ell ^1_{\hat n}$ denotes
the projection
   \[ (\beta _j)_{j\not= n} \mapsto (\cdots , 0, \beta _{-N}, \cdots ,
      \beta _N, 0, \cdots ) .
   \]
To show that $T^n : {\mathcal W} \rightarrow {\mathcal L}(\ell ^1_{\hat n})$
is weakly analytic it suffices to show that on any disk 
${\mathcal D}_{\varphi,h} := \{ \varphi + z h\,|\, z \in {\mathbb C}, |z| < 1\} $ with 
closure $\overline{{\mathcal D}_{\varphi,h}}\subseteq {\mathcal W}$, $\Pi_N \circ T^n$ converges 
in ${\mathcal L}(\ell^1_{\hat n})$ to $T^n$ locally uniformly in ${\mathcal D}_{\varphi,h}$ as 
$N \rightarrow \infty $. Indeed, if this is the case it follows from the Weierstrass theorem
that $T^n \big\arrowvert _{{\mathcal D}_{\varphi,h}} : {\mathcal D}_{\varphi,h}\rightarrow 
{\mathcal L}(\ell ^1_{\hat n})$ is analytic,
establishing in this way that $T^n$ is weakly analytic. To see that $\Pi_N \circ T^n$ converges 
locally uniformly on ${\mathcal D}_{\varphi,h}$ to $T^n$ as $N \rightarrow \infty $, observe that 
$\overline {\mathcal D}_{\varphi,h}$ is a compact subset of ${\mathcal W}$. The claimed
convergence thus follows from the estimates of Lemma~\ref{Lemma 7.1} and Lemma~\ref{Lemma 7.3}
(cf. the proof of Proposition \ref{Proposition 4.4}).

\noindent (iii) Arguing in the same way as above one shows that 
$T^\ast : {\mathcal W} \rightarrow {\mathcal L}(\ell ^1)$ is analytic.
\finishproof

\medskip

\noindent{\em Proof of Proposition~\ref{Proposition 4.3}.} 
Take $n\in\Z$. Let us begin by showing that 
$\beta^n = (T^n)^{-1} b^n$, $\beta^n : {\mathcal W} \rightarrow \ell^1_{\hat n}$, is analytic. 
As by Lemma~\ref{Lemma 7.4} (ii) $T^n : {\mathcal W}\rightarrow {\mathcal L}(\ell ^1_{\hat n})$
is analytic and by Proposition \ref{Proposition 4.2} (ii),  $T^n(\varphi)\in{\mathcal L}(\ell^1_{\hat n})$ is
a linear isomorphism  for any $\varphi\in\W$, it follows that
$(T^n)^{-1} : {\mathcal W} \mapsto {\mathcal L}(\ell ^1_{\hat n}), \varphi\mapsto (T^n(\varphi ))^{-1}$
is analytic as well. This combined with the analyticity of $b^n$,  
established in Lemma~\ref{Lemma 7.4} (i),  implies that $\beta^n$ is analytic.

Let us now turn towards $\zeta _n$. As the cases $|n| > N_0$ and $|n| \leq N_0$
are proved in the same way we consider only $|n| > N_0$. Then $\zeta _n$ is given by
   \[ \zeta _n : {\mathbb C} \times {\mathcal W} \rightarrow {\mathbb C}, \
      (\lambda , \varphi ) \mapsto \frac{\dot \Delta (\lambda )}{\lambda -
      \dot \lambda _n} - \xi ^n_{\beta ^n}(\lambda ) .
   \]
As $|n| > N_0$, $\dot \lambda _n : {\mathcal W} \rightarrow {\mathbb C}$ is
analytic and so is ${\mathbb C} \times {\mathcal W} \rightarrow {\mathbb C},(\lambda , \varphi ) 
\mapsto\frac{\dot \Delta (\lambda )}{\lambda - \dot \lambda _n}$, as $\lambda -\dot \lambda _n$ is a factor 
in the product representation for $\dot \Delta (\lambda )$ (Proposition \ref{Proposition 2.3}).
Recall that for any $\beta \in \ell ^1_{\hat n}$ and $|n| > N_0$,
   \begin{eqnarray*} \xi ^n_\beta (\lambda ) &=& \sum _{|j| > N_0}
      \beta _j \frac{\dot \Delta (\lambda )}{(\lambda - \dot \lambda _j)(\lambda -
      \dot \lambda _n)}\\
               &+& \Big( \sum ^{2N_0}_{j=0} \beta _{j-N_0}\lambda ^j \Big)
      \frac{\dot \Delta (\lambda )}{(\lambda - \dot \lambda _n) \prod_{|k| \leq
      N_0}(\lambda - \dot \lambda _k)} .
   \end{eqnarray*}
As above one argues that 
$\frac{\dot \Delta (\lambda )}{(\lambda - \dot \lambda_j)(\lambda - \dot \lambda _n)}$ and 
$\frac{\dot \Delta (\lambda )}{(\lambda -\dot \lambda _n)\prod _{|k| \leq N_0}(\lambda - \dot \lambda _k)}$ 
are analytic on ${\mathbb C} \times {\mathcal W}$. In addition, they are locally bounded uniformly in $j$.
By \cite[Theorem A.3]{KP} it then follows that 
the mapping ${\mathbb C} \times {\mathcal W} \rightarrow \ell ^\infty _{\hat n} $,
that assigns to $(\lambda,w)\in \C\times\W$ the sequence
   \[
      \left(\!\Big(\frac{\lambda^{N_0+j}}{\prod\limits_{|k|\leq N_0}
      (\lambda-\dot\lambda_k)}\frac{\dot\Delta(\lambda)}{\lambda-\dot\lambda_n}\Big)
               _{|j| \leq N_0},
      \Big(\frac{\dot\Delta(\lambda)}{(\lambda-\dot\lambda_n)(\lambda-\dot\lambda_j)}\Big)
                        _{|j|>N_0, j\ne n}\!\!\right)\!\in\!\ell^\infty_{\hat n}
   \]
is analytic. This combined with the result above, saying that 
$\beta^n: {\mathcal W} \rightarrow \ell ^1_{\hat n}$ is analytic it follows that
   \[ \xi ^n_{\beta ^n} : {\mathbb C} \times {\mathcal W} \rightarrow
      {\mathbb C} , \quad (\lambda , \varphi ) \mapsto \xi ^n_{\beta ^n
      (\varphi )} (\lambda , \varphi )
   \]
is analytic. Finally, by construction, the identity 
$\frac{1}{2\pi}\int_{A_m}\frac{\zeta_n(\lambda)}{\sqrt{{\mathcal R}(\lambda)}}\,d\lambda =
\delta _{nm}$ holds for any $n,m \in {\mathbb Z}$.
\finishproof

\medskip

\noindent To finish this section we prove results for the limiting behavior of
$b^n, T^n$, and $\beta ^n$ as $|n| \rightarrow \infty $. To this end introduce
the maps $\hat b ^n : {\mathcal W} \rightarrow \ell ^1$, 
$\hat T^n : {\mathcal W} \rightarrow {\mathcal L}(\ell ^1)$, and 
$\hat \beta ^n : {\mathcal W} \rightarrow \ell ^1$,
   \begin{align*} \hat b^n_m &:= \begin{cases} b^n_m , &m \not=
                     n \\ 0 , & m = n \end{cases} \qquad \quad
                     \hat \beta^n_j = \begin{cases} \beta^n_j , & 
                       j\not= n \\ 0 , & j = n \end{cases} \\
                  \hat T^n_{mj} &:= \begin{cases} T^n_{mj} , &
                     j, m \in {\mathbb Z} \backslash \{ n\} \\
                     0 , &(j,m) \in \Big(\big({\mathbb Z} \backslash\{ n\}\big)\times \{ n\}\Big)
                                       \cup \Big(\{ n \} \times \big({\mathbb Z}\backslash \{ n \}\big)\Big)\\ 
                     1 , & (j,m) = (n,n) .
                     \end{cases}
   \end{align*}

\begin{Lemma}\label{Lemma 7.5} 
For any given $\varphi\in\W$,
\begin{itemize}
\item[$(i)$] $\lim\limits_{|n|\to\infty }\|\hat b^n-b^*\| _{\ell ^1} = 0$;
\item[$(ii)$] $\lim\limits_{|n|\to\infty }\|\hat T^n - T^\ast\|_{{\mathcal L}(\ell^1)} = 0$;
\item[$(iii)$] $\lim\limits_{|n|\to\infty}$ $\|\hat\beta^n-\beta^*\|_{\ell ^1}=0$.
\end{itemize}
\end{Lemma}

\noindent{\em Proof.} Take $\varphi\in\W$. For simplicity of the notation we drop the argument
$\varphi$ in the quantities below.

\noindent (i) As $b^*_m=0$ for $|m|\ge N_0$ (see \eqref{eq:zero_form1}) we 
see from the definition of ${\hat b}^n_m$ that for $|n|\ge N_0$,
\[
{\hat b}^n_n-b^*_n=0\,.
\]
Moreover, by Lemma~\ref{Lemma 7.2},
   \[ |\hat b^n_m-b^*_m| \leq \begin{cases} C\,\frac{|\dot \lambda _m - \tau _m| +
      |\gamma _m|}{|m - n|} , &|m| > N_0,\, m \not= n \\ C\,
      \frac{1}{|m - n|} , &|m| \leq N_0,\, m \not= n . \end{cases}
   \]
By Proposition~\ref{Proposition 2.2}, 
$|\dot \lambda _m - \tau _m| + |\gamma_m| = \ell ^2(m)$. 
Therefore, for $|n|\ge N_0$,  $|\hat b^n_m-b^*_m| \leq C \frac{a_m}{|m-n|}$
where $(a_m)_{m \in {\mathbb Z}} \in \ell ^2$. 
Thus, for $|n|\ge N_0$,
   \begin{align} \sum _{m \in {\mathbb Z}} |{\hat b}^n_m-b^*_m| &\le C\!\!\!\!\!\sum
                     _{|m - n| \leq | \frac{n}{2}|, m\ne n}
                     \frac{a_m}{|m - n|} + C\!\!\sum _{|m - n| > |\frac{n}{2}|} 
                     \frac{a_m}{|m - n|}\nonumber\\
                  &\leq C \Big( \sum _{|m - n| \leq | \frac{n}{2}|}
                     a_m^2 \Big) ^{1/2} \Big( \sum _{k\not= 0} \frac{1}{k^2}
                     \Big) ^{1/2}\nonumber\\
                  &+ C \Big( \sum _{m \in {\mathbb Z}} a_m^2 \Big)^{1/2}
                     \Big( \sum _{k \not= 0} \frac{1}{k ^{4/3}} \Big) ^{1/2}
                     \Big(\frac{2}{|n|}\Big)^{1/3}
   \end{align}
implying that $\lim\limits_{|n|\to\infty}\|{\hat b}^n-b^*\|_{\ell^1}=0$.

\medskip

\noindent (ii) Using the same arguments as in the prove of (i) one concludes
from Lemma~\ref{Lemma 7.3} that the claimed convergence holds.

\medskip

\noindent (iii) It follows from Proposition \ref{Proposition 4.2} and 
Proposition \ref{Proposition 4.4} that $\hat T^n$  and $T^*$ are linear isomorphisms 
in $\ell^1$. Hence (ii) implies that for any $\varphi\in\W$,
$(\hat T^n)^{-1} \to (T^\ast )^{-1}$ in ${\mathcal L}(\ell ^1)$
as $|n|\to\infty$. This together with (i) and $\hat \beta ^n = (\hat T^n)^{-1}{\hat b}^n$ imply
that for any $\varphi\in\W$,
\[
{\hat\beta}^n=(\hat T^n)^{-1}{\hat b}^n\to(\hat T^*)^{-1} b^*=\beta^*
\]
in $\ell^1$ as $|n|\to\infty$.
\finishproof


\section{Estimates of the zeros}\label{Estimates of the zeros}

In this section we prove that the zeros of the analytic function $\zeta _n :
{\mathbb C} \times {\mathcal W} \rightarrow {\mathbb C}$, introduced in
Section~\ref{3. Outline of proof of Theorem 1.2}, satisfy the properties stated
in Theorem~\ref{Theorem 1.2}. The ansatz we have chosen for
the $\zeta _n$ is well suited to obtain these claimed estimates. Recall from \eqref{eq:zeta_n} that
for any $|n| > N_0$
   \[ \zeta _n(\lambda ) \equiv \zeta _n(\lambda , \varphi ) = \big(1 - \eta _n(\lambda)\big) 
        \frac{\dot \Delta (\lambda )}{\lambda - \dot \lambda _n}
   \]
where
   \[ \eta _n(\lambda ):= \sum _{|j| > N_0, j\ne n} \frac{\beta
      ^n_j}{\lambda - \dot \lambda _j} + \frac{p_n(\lambda )}{\prod _{|j|
      \leq N_0}(\lambda - \dot \lambda _n)}
   \]
and $p_n(\lambda ):=p^n_{\beta ^n}(\lambda , \varphi )$ is the polynomial
introduced in Section~\ref{3. Outline of proof of Theorem 1.2} with $\beta $
given by $\beta ^n \equiv \beta ^n(\varphi )$ of Proposition~\ref{Proposition 4.3}. 
Similarly, for $|n|\le N_0$ we define
 \[ \zeta _n(\lambda ) \equiv \zeta _n(\lambda , \varphi ) = \big(1 - \eta _n(\lambda)\big) 
        \frac{\dot \Delta (\lambda )}{\lambda - \dot \lambda_{N_0}}
   \]
where
   \[ \eta _n(\lambda ):= \sum _{|j| > N_0} \frac{\beta
      ^n_j}{\lambda - \dot \lambda _j} + \frac{p_n(\lambda )}{\prod _{|j|
      \leq N_0, j\ne N_0}(\lambda - \dot \lambda _j)}\,.
   \]
First note that for any $n \in {\mathbb Z}$, $\frac{\dot \Delta (\lambda )}{\lambda - \dot \lambda _n}$
is an entire function and by the argument principle one has in view of the choice of $\W$ and
Proposition~\ref{Proposition 2.1}  that for any $\varphi \in {\mathcal W}$
   \begin{equation}
   \label{8.1} \frac{1}{2\pi i} \int _{\Gamma _m} \partial _\lambda \log
               \Big( \frac{\dot \Delta (\lambda )}{\lambda - \dot \lambda _n}
               \Big)\, d \lambda = 1 - \delta _{nm} \quad \forall |m| > N_0
   \end{equation}
and
   \begin{equation}
   \label{8.2} \dot \lambda _m = \frac{1}{2\pi i} \int _{\Gamma _m} \lambda\,
               \partial _\lambda \log \Big( \frac{\dot \Delta (\lambda )}
               {\lambda - \dot \lambda _n} \Big)\, d\lambda \quad \forall
               |m| > N_0, m\ne n
   \end{equation}
whereas for any $N \geq N_0$
   \begin{equation}
   \label{8.3} \frac{1}{2\pi i} \int _{\partial D_0(N\pi + \frac{\pi }{4})}
               \partial _\lambda \log \Big( \frac{\dot \Delta (\lambda )}
               {\lambda - \dot \lambda _n}\Big)\, d\lambda = \begin{cases} 2N,
               &|n| \leq N \\ 2N + 1 , &|n| > N .
               \end{cases}
   \end{equation}
Viewing $\zeta _n(\lambda )$ as a perturbation of $\frac{\dot \Delta (\lambda )}
{\lambda - \dot \lambda _n}$ we want to argue in a similar fashion for $\zeta
_n(\lambda )$. First we need to establish some auxiliary estimates.

\begin{Lemma}\label{Lemma 8.1} 
For any $\varphi \in {\mathcal W}$, $\beta \in \ell ^1$,  $N > N_0\ge 1$, 
$n \in {\mathbb Z}$, and $|m| > 2N$, one has
\begin{itemize}
\item[$(i)$] $\sup\limits_{\lambda \in \Gamma _m}\Big(\sum _{|j| > N_0} \frac{|\beta _j|}
{|\lambda - \dot \lambda _j|}\Big) \leq C \Big( |\beta _m| + \sum\limits_{|j| > N, j\ne m}
\frac{|\beta _j|}{|j - m|} + \| \beta \| _{\ell ^1} \frac{1}{m}\Big)$;
\item[$(ii)$] $\sup\limits_{\lambda \in \Gamma _m}\Big(\frac{|p^n_\beta (\lambda )|}{\prod
_{|j| \leq N_0} |\lambda - \dot \lambda _j|}\Big) \leq C \| \beta \| _{\ell ^1}
\frac{1}{m}$, \quad $|n| > N_0$;
\item[$(iii)$] $\sup\limits_{\lambda \in \Gamma _m} \Big(\frac{|p^n_\beta (\lambda )|}
{\prod _{|j| \leq N_0, j\ne N_0}|\lambda - \dot \lambda _j|}\Big)
\leq C \| \beta \| _{\ell ^1} \frac{1}{m}$, \quad $|n| \leq N_0$\,,
\end{itemize}
where the constant $C > 0$ can be chosen uniformly  in $N>N_0$, $n\in{\mathbb Z}$, 
$|m|>2N$, and $\varphi\in{\mathcal W}$.
\end{Lemma}

\noindent{\em Proof. } For any $|m| > 2N$ and $\lambda \in \Gamma _m$ one has
by the choice of $\W$ and Proposition~\ref{Proposition 2.1}, $1 / |\lambda - \dot \lambda _m|\leq \pi / 12$ as 
$\mbox{dist}(\dot \lambda _m, \Gamma _m) \ge\frac{\pi}{4}-\frac{\pi}{6}$. In view
of the choice of $\W$ and Proposition~\ref{Proposition 2.1} it then follows that for $|m| > 2N$
and $\lambda \in \Gamma _m$
   \begin{align*} \sum _{|j| > N_0} \frac{|\beta _j|}{|\lambda - \dot \lambda _j|}
                     &\leq C \Big( |\beta _m| + \sum_{|j|> N, j\ne m}
                      \frac{|\beta _j|}{|j-m|}\Big) + \sum _{N_0 < |j| \leq N}
                     |\beta _j| \frac{1}{|\lambda - \dot \lambda _j|} \\
                  &\leq C \Big( |\beta _m| + \sum _{|j|> N, j\ne m} 
                        \frac{|\beta _j|}{|j-m|}\Big) + C \| \beta \|_{\ell ^1}\frac{1}{m}
   \end{align*}
where $C > 0$ can be chosen uniformly in $N>N_0$, $|m|>2 N$, and $\varphi \in {\mathcal W}$. 
Towards (ii) note that as $|\lambda|\ge 1$ one has
   \[ \frac{|p^n_\beta (\lambda )|}{\prod _{|j| \leq N_0} |\lambda - \dot
      \lambda _j|} \leq \frac{1}{|\lambda |} \frac{\sum _{|j| \leq N_0} |\beta
      _j|}{\prod _{|j| \leq N_0}\Big|1 - \frac{\dot \lambda _j}{\lambda }\Big|}
   \]
for any $|n|>N_0$ and $\varphi\in\W$.   
In addition $|\dot \lambda _j| < (N_0 + \frac{1}{4}) \pi $ for any 
$|j| \le N_0$ and $|\lambda | \geq m\pi - \pi / 4$. Hence,
$\prod _{|j| \leq N_0}\Big|1 - \frac{\dot \lambda _j}{\lambda }\Big| \geq C$ 
where the constant $C > 0$ can be chosen uniformly in $\varphi \in {\mathcal W}$. 
This implies that
   \[ \sup\limits_{\lambda \in \Gamma _m}\Big(\frac{|p^n_\beta (\lambda )|}{\prod _{|j|
      \leq N_0 }|\lambda - \dot \lambda _j|}\Big) \leq C\| \beta \|_{\ell^1}\frac{1}{m}
   \]
uniformly in $|n|>N_0$, $N>N_0$, $|m|>2 N$, and $\varphi\in\W$.

\noindent Item (iii) is proved in a similar fashion.
\finishproof

Next we want to estimate $\eta _n(\lambda )$ on $\partial D_0(r_{2N})$ where
for any $m \in {\mathbb Z}$, $r_m := m\pi + \pi /4$.

\begin{Lemma}\label{Lemma 8.2} 
For any $\varphi \in {\mathcal W}$, $\beta\in {\ell }^1$, $N > N_0\ge 1$, and
$n \in {\mathbb Z}$, one has
\begin{itemize}
\item[$(i)$] $\sup\limits_{\lambda \in \partial D_0(r_{2N})}\Big(\sum _{|j| > N_0}
\frac{|\beta _j|} {|\lambda - \dot \lambda _j|}\Big) \le C \frac{\|\beta \|_{\ell^1}}{N} +  
C \sum _{|j| > N} | \beta _j|$;
\item[$(ii)$] $\sup\limits_{\lambda \in \partial D_0(r_{2N})}\Big(\frac{|p^n_\beta
(\lambda )|}{\prod _{|j| \leq N_0} |\lambda - \dot \lambda _j|}\Big) \leq C \| \beta
\| _{\ell ^1} \frac{1}{N}$, \quad $|n| > N_0$;
\item[$(iii)$] $\sup\limits_{\lambda \in \partial D_0(r_{2N})}\Big(\frac{|p^n_\beta
(\lambda )|}{\prod\limits_{|j| \leq N_0, j\ne N_0} |\lambda - \dot \lambda_j|}\Big) 
\leq C \| \beta \| _{\ell ^1} \frac{1}{N}$, \quad $|n| \leq N_0$,
\end{itemize}
where $C > 0$ can be chosen uniformly in $n\in\Z$, $N>N_0$, and $\varphi \in {\mathcal W}$.
\end{Lemma}

\noindent{\em Proof.} To prove item (i) we split the sum $\sum _{|j| > N_0}$ into two parts: 
$\sum _{N_0 < |j| \leq N}$ and $\sum _{|j| > N}$. Clearly, for any $\lambda\in \partial D_0(r_{2N})$,
\[
\sum _{N_0 < |j| \leq N} \frac{|\beta _j|}{|\lambda -\dot \lambda _j|} \leq 
C \| \beta \|_{\ell ^1}/N
\]
and
   \[ \sum _{|j| > N} \frac{|\beta _j|}{|\lambda - \dot \lambda _j|} \leq
      \Big( \sum _{|j| > N} |\beta _j|^2 \Big)^{1/2} \Big( \sum _{|j| > N}
      \frac{1}{|\lambda - \dot \lambda _j|^2} \Big) ^{1/2} \leq C \sum _{|j|
      > N} |\beta _j|
   \]
where $C > 0$ can again be chosen independently of $n \in {\mathbb Z}$, $N>N_0$, and $\varphi \in {\mathcal W}$. 
The estimates (ii) and (iii) are proved in the same
way as items (ii) respectively (iii) of Lemma~\ref{Lemma 8.1}.
\finishproof

Lemma~\ref{Lemma 8.1} and Lemma~\ref{Lemma 8.2} can be used to localize the
zeros of $\zeta _n(\cdot , \varphi )$ for any $n \in {\mathbb Z}$ and $\varphi\in {\mathcal W}$. 
In fact, let $\psi\in L^2_\bullet$ be the potential appearing in the construction of $\W$ in 
Section \ref{2. Preliminaries}. Using that $\beta^n(\psi)\to\beta^*(\psi)$ in $\ell^1$ as $|n|\to\infty$
(Lemma \ref{Lemma 7.5}) and the fact that $\beta^n : \W\to\ell^1_{\hat n}$ and 
$\beta^* : \W\to\ell^1$ are analytic (Corollary \ref{coro:beta_analytic}), we conclude that
for any $\varepsilon>0$ there exists $N_1\ge 1$ and an open neighborhood $\W_1\subseteq\W$ 
of $\psi$ in $L^2_\bullet$ such that for any $n\in\Z$ and for any $\varphi\in\W_1$,
\begin{equation}\label{eq:beta_bounded}
\sum_{|m|\ge N_1}|{\hat\beta}^n_m|<\varepsilon\,.
\end{equation}
Combining this with Lemma \ref{Lemma 8.1} and Lemma \ref{Lemma 8.2}, shrinking the neighborhood $\W$
and taking $N_1\ge N_0\ge 1$ larger if necessary we see that 
for any $N \geq N_1$, $n\in{\mathbb Z}$, $|m| > 2N$, and $\varphi\in\W$,
   \begin{equation}
   \label{8.4} \sup _{\lambda \in \Gamma _m} |\eta _n(\lambda )| \leq 1/2
   \end{equation}
and
   \begin{equation}
   \label{8.5} \sup _{\lambda \in \partial D _0(r_{2N})} |\eta _n(\lambda )|
               \leq 1/2\,.
   \end{equation}
It then follows that for any $N\ge N_1$, $\lambda \in \Gamma _m$, $|m| > 2N$, as well as 
for any $\lambda\in \partial D_0(r_{2N})$,
\[
\left\{
\begin{array}{l}
\,\,\Big|\frac{\dot\Delta (\lambda )}{\lambda - \dot
                     \lambda _n} - \zeta_n(\lambda )\Big| \leq
                     \Big| \eta _n(\lambda ) \frac{\dot \Delta (\lambda )}
                     {\lambda - \dot \lambda _n} \Big| \leq \frac{1}
                     {2} \Big| \frac{\dot \Delta (\lambda )}{\lambda -
                     \dot \lambda _n}\Big|,\,\,\,\,\,|n| >N_0 \\
\Big|\frac{\dot \Delta (\lambda )}{\lambda - \dot
                     \lambda _{N_0}} - \zeta_n(\lambda )\Big| \leq
                     \Big| \eta _n(\lambda ) \frac{\dot \Delta (\lambda )}
                     {\lambda - \dot \lambda _{N_0}} \Big| \leq \frac{1}
                     {2} \Big| \frac{\dot \Delta (\lambda )}{\lambda -
                     \dot \lambda _{N_0}}\Big|,\,\,\,\,\,|n| \leq
                     N_0\,.
\end{array}
\right.
\]
Hence by Rouch\'e's theorem and formulas \eqref{8.1} and \eqref{8.3} one
has for any $N>N_1$ and $|m| > 2N$,
   \[ \frac{1}{2\pi i} \int _{\Gamma _m} \partial _\lambda \big( \log \zeta _n
      (\lambda ) \big)\,d\lambda = 1 - \delta _{nm}
   \]
   \[ \frac{1}{2\pi i} \int _{\partial D_0(r_{2N})} \partial _\lambda \big(
      \log \zeta _n(\lambda ) \big)\,d\lambda = \begin{cases} 4N , &
      |n| \leq 2N \\ 4N + 1 , &|n| > 2N . \end{cases}
   \]
\begin{Rem}
As these relations hold for {\em any} $N>N_1$ we also see that for any $\varphi\in\W$  and for any $N>N_1$
there are no zeros of $\zeta_n(\cdot,\varphi)$ outside of the union of the sets $D_m(\pi/4)$, $|m|>2N$, and
$D_0(r_{2N})$.
\end{Rem}
For any $|m| > 2N$ we denote the zero of $\zeta _n(\lambda )$ inside
$\Gamma _m$ by $\sigma ^n_m$. By the argument principle one has (cf. \eqref{8.2})
\begin{align*}
 \begin{split}
        \sigma ^n_m &= \frac{1}{2\pi i} \int _{\Gamma _m} \lambda\,
                     \partial _\lambda \big( \log \zeta _n(\lambda ) \big)\,
                     d \lambda \\
                    &= \dot \lambda _m + \frac{1}{2\pi i} \int _{\Gamma _m}
                  \lambda\,\partial _\lambda \big( \log (1 - \eta _n(\lambda ))
                  \big)\, d\lambda\,.
 \end{split}
\end{align*}
Integration by parts leads to
\begin{equation}\label{8.6} 
\sigma ^n_m = \dot \lambda _m - \frac{1}
{2\pi i} \int _{\Gamma _m} \log \big( 1 - \eta _n(\lambda )\big)\,d\lambda\,.
\end{equation}
Using \eqref{8.4} and Lemma~\ref{Lemma 8.1} one sees that for $|m|>2 N$,
   \begin{align*} &\sup _{\lambda \in \Gamma _m} \Big\arrowvert \log
                     \big( 1 - \eta _n (\lambda ) \big) \Big\arrowvert \leq 2
                     \sup _{\lambda \in \Gamma _m} \big\arrowvert\eta_n
                     (\lambda ) \big\arrowvert\le \\
                  &\leq C_1 \Big( |\beta ^n_m| + \sum _{|j|> N_1, j\ne m, n} 
                     \frac{|\beta ^n_j|}{|j -  m|} + \| \beta ^n \|_{\ell_{\hat n}^1}/m\Big)
   \end{align*}
where $C_1>0$ is independent of $n\in\Z$, $\varphi\in\W$, and $|m|>2 N$.
By using the Cauchy-Schwartz inequality and then changing the order of summation in the double sum we get 
\begin{align*}
\sum_{|m|>2N}\Big(\sum_{j\ne m,n}\frac{|\beta^n_j|}{|j-m|}\Big)^2\le
\sum_{|m|>2N}\|\beta^n\|_{\ell^1_{\hat n}}\sum_{j\ne m,n}\frac{|\beta^n_j|}{|j-m|^2}=\\
=\|\beta^n\|_{\ell^1_{\hat n}}\sum_{j\ne n}|\beta^n_j|\sum_{|m|>2N, m\ne j}\frac{1}{|j-m|^2}
\le C_2 \|\beta^n\|_{\ell^1_{\hat n}}^2
\end{align*}
where $C_2=2\sum_{k\ge 1}\frac{1}{k^2}$. Hence, we get in view of \eqref{eq:beta_bounded}
that there exists $C>0$ such that for any $n\in\Z$ and for any $\varphi\in\W$
\[
\sum _{|m| > 2N} |\sigma ^n_m - \dot \lambda _m|^2 \leq C\,.
\] 
Note that as $\zeta _n : {\mathbb C} \times {\mathcal W} \rightarrow {\mathbb C}$ is analytic by 
Proposition~\ref{Proposition 4.3}, the identity \eqref{8.6} also shows that 
$\sigma ^n_m : {\mathcal W} \rightarrow {\mathbb C}$ is analytic for any $|m| > 2N$, $m\ne n$. 
By denoting $2N$ again by $N$, we get
\begin{Prop}\label{Proposition 8.3} 
For any $\psi\in L^2_\bullet$ there exists an open neighborhood $\W$ of $\psi$ in $L^2_\bullet$
obtained from Section \ref{2. Preliminaries} after shrinking if necessary and $N\geq N_0$ so that for any $n\in{\mathbb Z}$ and for any 
$\varphi \in {\mathcal W}$,  the entire function $\zeta_n(\lambda)$ has precisely $2N + 1$  [$2N$] zeros 
inside $D_0(r_N)$ if $|n| > N$  [$|n| \leq N$]. For any $|m| > N$, $m\ne n$, $\zeta _n(\lambda)$ 
has precisely one zero,  denoted by $\sigma^n_m \equiv \sigma ^n_m(\varphi)$, in $D_m(\pi/4)$. 
There are no other zeros of $\zeta_n(\lambda)$ in $\C$. 
Moreover, $\sigma^n_m={\dot\lambda}_m+\ell^2(m)$, $|m|>N$, uniformly in $n \in {\mathbb Z}$ and 
uniformly in ${\mathcal W}$, and for any $|m| > N$, 
$\sigma ^n_m : {\mathcal W} \rightarrow {\mathbb C}$ is analytic.
\end{Prop}
Proposition~\ref{Proposition 8.3} implies that $\zeta _n(\lambda )\equiv \zeta _n(\lambda,\varphi )$ has 
in fact a product representation. List the roots of $\zeta _n(\lambda,\varphi )$ inside $D_0(r_N)$ in
lexicographic order and with their multiplicities, $\sigma ^n_m$, $|m| \le N$, $m\ne n$.
\begin{Coro}\label{Corollary 8.4} 
For any $n \in {\mathbb Z}$ and $\varphi \in {\mathcal W}$,
   \[ \zeta _n(\lambda ) = - \frac{2}{\pi _n} \prod _{j \not= n} \frac{\sigma ^n
      _j - \lambda }{\pi _j} .
   \]
\end{Coro}
\noindent{\em Proof.} Take $\varphi\in\W$. As the cases $|n| > N_0$ and $|n| \leq N_0$ are proved in
the same way let us consider the case $|n| > N_0$. By \cite[Theorem 2.2]{GKP},
$\dot \Delta (\lambda )$ is an entire function of order 1 and so is $\frac{\dot
\Delta (\lambda )}{\lambda - \dot \lambda _n}$. By Proposition~\ref{Proposition 4.3}, 
$\zeta _n(\lambda )$ is an entire function and by Lemma~\ref{Lemma 8.2},
   \[ \sup _{\lambda \in \partial D_0(r_{2N})} |1 - \eta _n(\lambda )| = 1+o(1)
      \mbox { as } N \rightarrow \infty  .
   \]
It then follows that 
$\zeta _n(\lambda ) = \big(1 - \eta _n(\lambda )\big) \frac{\dot\Delta (\lambda )}{\lambda - \dot \lambda _n}$ 
is an entire function of order $1$. Moreover, by Proposition \ref{Proposition 2.2} and 
Proposition \ref{Proposition 8.3} the exponent of convergence of the zeros of $\zeta_n(\lambda)$ is
equal to $1$ and the series $\sum_{|j|>N_0}\frac{1}{|\sigma^n_j|}$ diverges.
This implies that the genus of $\zeta_n(\lambda)$ is equal to $1$.
By Hadamard's factorization theorem
\begin{equation}\label{eq:hadamard}
\zeta _n(\lambda ) = \lambda ^{\nu_n}\,e^{\tilde a_n \lambda + \tilde b_n}\,
\prod_{\sigma ^n_k \not= 0} E\Big( \frac{\lambda }{\sigma ^n_k},1\Big)
\end{equation}
where $\nu _n$ is the order of vanishing of $\zeta _n(\lambda )$ at $\lambda= 0$, 
$\tilde a_n, \tilde b_n \in {\mathbb C}$ are constants independent of $\lambda\in\C$,
and $E(z,1)$ is the canonical factor $E(z,1):=(1 - z )e^z $.
For $|m| > n$ we pair the factors 
$E\Big(\frac{\lambda }{\sigma ^n_{-m}},1\Big) \cdot E\Big(\frac{\lambda }{\sigma ^n_m},1\Big)$ and
conclude from \eqref{eq:hadamard}, Proposition \ref{Proposition 2.2}, and Proposition \ref{Proposition 8.3},
that $\zeta_n(\lambda )$ has a product representation of the form
   \[ e^{a_n \lambda + b_n} \prod _{k\ne n}
      \frac{\sigma ^n_k - \lambda }{\pi _k} = e^{a_n \lambda + b_n}
      \lim_{K\to\infty}\prod_{|k|\le K, k\ne n}\frac{\sigma^n_k-\lambda}{\pi_k}\,.
   \]
On the other hand, by Proposition~\ref{Proposition 2.3}, $\frac{\dot \Delta
(\lambda )}{\lambda - \dot \lambda _n} = - \frac{2}{\pi _n} \prod _{k \not=n} 
\frac{\dot \lambda _k - \lambda }{\pi _k}$ and by \cite[Lemma C.5]{GKP},
on the circles $|\lambda | = r_{2N}$,
   \[ \frac{\prod _{k \not= n} \frac{\sigma ^n_k - \lambda}{\pi _k}}
      {\prod _{k \not= n} \frac{\dot \lambda _k - \lambda}{\pi _k}}
      = 1 + o(1) \quad \mbox{ as } N \rightarrow \infty \,.
   \]
By Lemma \ref{Lemma 8.2} and Proposition \ref{Proposition 8.3}, on the circle $|\lambda | = r_{2N}$,
   \[ \frac{e^{a_n\lambda + b_n}}{- \frac{2}{\pi _n}} \frac{\prod _{k \not= n}
      \frac{\sigma ^n_k - \lambda}{\pi _k}}{\prod _{k \not= n} \frac{\dot
      \lambda _k - \lambda }{\pi _k}} = 1 + o(1) \quad \mbox{ as } N \rightarrow
      \infty .
   \]
It then follows that $a_n = 0$ and $e^{b_n} = - \frac{2}{\pi _n}$, yielding
the claimed formula
$\zeta _n(\lambda ) = - \frac{2}{\pi _n} \prod _{k \not= n} 
\frac{\sigma ^n_k - \lambda }{\pi _k}$.
\finishproof

The refined asymptotics of the zeros $(\sigma ^n_m)_{m \not= n}$ of $\zeta _n$
stated below are proved in the same way as in \cite[Lemma 14.12]{GKP} and hence
we omit its proof.
\begin{Lemma}\label{Lemma 8.5} 
There exist $N\geq N_0$ so that
\[ 
\sigma ^n_m = \tau _m + \gamma ^2_m \ell ^2_m \quad\quad\forall |m| > N
\]
uniformly in $n\in\Z$ and locally uniformly in $\W$.
\end{Lemma}

\noindent Finally we prove the following
\begin{Lemma}\label{Lemma 8.6} 
For any $n \in {\mathbb Z}$ and $\varphi \in {\mathcal W}\subseteq L^2_\bullet$,
the entire function $\zeta _n(\lambda)$ vanishes at 
$\lambda\in Z_\varphi\setminus\{\lambda^\pm_n(\varphi)\}$. If $\lambda ^\pm_n \in Z_\varphi$
then $\zeta_n(\lambda)$ does not vanish at $\lambda ^\pm_n$.
\end{Lemma}
\noindent{\em Proof.} Take $\varphi\in\W$. To see that $\zeta _n(\lambda,\varphi )$ vanishes on
$Z_\varphi \backslash \{ \lambda ^\pm _n(\varphi )\} $ one argues as in the proof of Lemma~\ref{Lemma 5.2}. 
If $\lambda^\pm_n \in Z_\varphi$, then $\lambda^\pm_n$ is a zero of order two of ${\mathcal R}(\lambda)$
and hence $\frac{\zeta _n(\lambda )}{\sqrt[c]{{\mathcal R}(\lambda)}}\,d\lambda$ has a pole of order
$\le 1$ at $\lambda^\pm_n$. As 
$\int_{\Gamma_n}\frac{\zeta _n(\lambda )}{\sqrt[c]{{\mathcal R}(\lambda)}}\,d\lambda=2\pi$ we conclude that
$\zeta_n(\lambda^\pm_n)\ne 0$.
\finishproof


\end{document}